\documentclass[review]{elsarticle}
\usepackage{lineno,hyperref}
\usepackage{graphicx}
\modulolinenumbers[5]
\usepackage{booktabs}
\usepackage{multirow}
\usepackage{siunitx}
\usepackage{multirow}
\usepackage{graphicx,amsmath,amssymb}
\usepackage{amsmath}
\usepackage{epstopdf} 
\usepackage{xcolor}
\usepackage{pbox}
\usepackage{subfigure}
\usepackage{lineno,hyperref}
\usepackage{booktabs}
\usepackage{multirow}
\usepackage{siunitx}
\usepackage{mathptmx}      
\usepackage{algorithm}
\usepackage{algorithmic}
\usepackage{multirow}
\usepackage{multirow}
\usepackage{siunitx}
\usepackage[para]{footmisc}
\journal{Journal of \LaTeX\ Templates}

\begin{document}

\begin{frontmatter}

\title{Second order difference approximation for a class of Riesz space fractional
advection-dispersion equations with delay}

\author{M. Saedshoar Heris $^1$}
\author[M. Javidi ]{M. Javidi\corref{mycorrespondingauthor}}
\cortext[mycorrespondingauthor]{Corresponding author}
\ead{mo$_{-}$javidi@tabrizu.ac.ir}

\address{$^1$Department of Applied Mathematics, University of Tabriz, Tabriz, Iran}

\begin{abstract}
In this paper, we propose numerical scheme for the Riesz space fractional advection-dispersion equations with delay (RFADED). Firstly, analytical solution for RFADED in terms of the functions of Mittag-Leffler type is derived. Secondly, the fractional backward differential formulas (FBDF) method and shifted Gr\"{u}nwald method are introduced to the Riesz space fractional derivatives. Next, stability and convergency of these methods have been proved. Thirdly, Crank-Nicolson scheme for solving this problem is proposed. We prove that the scheme is conditionally stable and convergent with the order accuracy of ${\rm O}({\kappa ^2} + {h^2})$. Finally, some numerical results are given to demonstrate that presented method is a computationally efficient and accurate method for solving RFADED.
\end{abstract}


\end{frontmatter}


\section{Introduction}
In the past few decades, fractional calculus has attracted increasing interest due to its applications in science, for examples, physics, engineering, biology, economics and finance \cite{bagley1991fractional, weaver1990vibration, marks1981differintegral, simo2016effects}.
The fractional-order models have proved to be more accurate than integer-order models for the memory and hereditary properties of the fractional-order derivative \cite{podlubny1999fractional,kilbas2006theory}. Fractional differential equations (FDEs) provide a powerful and flexible tool for modeling and describing the behavior of science such as real materials \cite{bagley1991fractional}, fluid dynamics \cite{marks1981differintegral}, finance, electromagnetic waves, electrochemical process \cite{weaver1990vibration}, biological systems, control theory \cite{wu2012theory} and so on.
Some numerical methods to solve the FDEs have been discussed in \cite{galeone2006multistep, garrappa2015trapezoidal, gorenflo1997fractional, lubich1986discretized,heris2019predictor,javidi2019predictor}. A delay differential equation is a differential equation in which the derivative of the function at any time depends on the solution at the previous time. The fractional differential equations with delay (FDDEs) have important applications in many fields, such as, chemistry, physics, biology and finance \cite{wang2011analysis}. Analytical solutions of FDDEs are complicated and hence the numerical solutions play an important role. Numerical methods for solving FDDEs are very important. For more details we refer the reader to a series of articles \cite{heris2017fractional,heris2017fbdf5,morgado2013analysis,vcermak2016stability,wang2013numerical}. Partial differential equations with delays appear in many fields such as in biology, medicine, population ecology, control systems and climate models \cite{wu2012theory}. Some numerical methods for solving delay partial differential equations are studied in \cite{zubik2000chebyshev, jackiewicz2006spectral,tanthanuch2012symmetry}. Numerical methods for solving fractional partial differential equations are investigated in \cite{tadjeran2006second,liu2015finite,ding2013numerical}.\\
In this paper, we consider a class of Riesz space fractional advection-dispersion equations (RFADEs)\cite{yang2010numerical} with delay:
\begin{equation}\label{eq4}
\dfrac{{\partial u(x,t)}}{{\partial t}} +{K_\gamma }\dfrac{{{\partial ^\gamma }u(x,t - \tau)}}{{\partial {{t}^\gamma }}}= {K_\alpha }\dfrac{{{\partial ^\alpha }u(x,t)}}{{\partial {{\left| x \right|}^\alpha }}} + {K_\beta }\dfrac{{{\partial ^\beta }u(x,t)}}{{\partial {{\left| x \right|}^\beta }}} + f(x,t),
\end{equation}
subject to the initial condition:
\begin{equation}\label{eq5}
\begin{array}{l}
u(x,t) = g(x,t),\,\,\, - \tau  \le t \le 0,\,\,\,0 \le x \le L,\\
u(0,t) = {\mu _1}(t),\,\,u(L,t) = {\mu _2}(t),\,\,0 \le t \le T,
\end{array}
\end{equation}
where $0 < \gamma  < 1,\,\,0 < \alpha  < 1,\,\,1 < \beta  \le 2$. $u$ is a solute concentration, ${K_\alpha}$ and ${K_\beta}$ are the dispersion coefficient and the average fluid velocity, respectively. Here, we take ${K_\alpha > 0 },\,\,\,{K_\beta \ge 0},\,\,\,{K_\gamma \ge 0 }$. Also $f$ is a continuous function on $[0,T] \times [0,L]$, $T>0, L>0$ and $g$ is a continuous function on $[ - \tau ,0] \times [0,L]$, $\tau>0$.
The Riesz space fractional operator on a finite domain [0, L] is defined as \cite{yang2010numerical}
\begin{equation}\label{eq6}
\dfrac{{{\partial ^\zeta }u(x,t)}}{{\partial {{\left| x \right|}^\zeta }}} =- {( - \Delta )^{\dfrac{\zeta }{2}}}u(x,t)= - {c_\zeta }[{}_0^{RL}D_x^\zeta u(x,t) + {}_x^{RL}D_L^\zeta u(x,t)],
\end{equation}
where
\begin{equation}\label{eq6666}
{c_\zeta } = \dfrac{1}{{2\cos (\frac{{\pi \zeta }}{2})}},\,\,\,\,0 < \zeta  \le 2,\,\,\zeta  \ne 1.
\end{equation}
Also, the left and right Riemann-Liouville derivatives of order $0<\zeta \le 2$, for a given function $u(x, t)$ on a finite interval [0, L], are defined as \cite{kilbas2006theory}:
\begin{equation}\label{eq40000}
\begin{array}{l}
{}_0^{RL}D_x^\zeta u(x,t) = \dfrac{1}{{\Gamma (1 - \zeta )}}\dfrac{\partial }{{\partial x}}\int\limits_0^x {{{(x - \eta )}^{ - \zeta }}u(\eta ,t)d\eta ,} \,\,0 < \zeta  < 1,\\
{}_x^{RL}D_L^\zeta u(x,t) = \dfrac{{ - 1}}{{\Gamma (1 - \zeta )}}\dfrac{\partial }{{\partial x}}\int\limits_x^L {{{(\eta  - x)}^{ - \zeta }}u(\eta ,t)d\eta ,} \,\,0 < \zeta  < 1,\\
{}_0^{RL}D_x^\zeta u(x,t) = \dfrac{1}{{\Gamma (2 - \zeta )}}\dfrac{{{\partial ^2}}}{{\partial {x^2}}}\int\limits_0^x {{{(x - \eta )}^{1 - \zeta }}u(\eta ,t)d\eta ,} \,\,1 < \zeta  \le 2,\\
{}_x^{RL}D_L^\zeta u(x,t) = \dfrac{1}{{\Gamma (2 - \zeta )}}\dfrac{{{\partial ^2}}}{{\partial {x^2}}}\int\limits_x^L {{{(\eta  - x)}^{1 - \zeta }}u(\eta ,t)d\eta ,} \,\,1 < \zeta  \le 2.
\end{array}
\end{equation}
Note that the time fractional derivative $\dfrac{{{\partial ^\gamma }}}{{\partial {{t}^\gamma }}}$ at the Eq (\ref{eq4}) is the left Riemann-Liouville derivatives.
Analytical and numerical solution of partial fractional differential equations (PFDEs) with delay are more complicated and less studied in the literature. The solutions of PFDEs are different from the ones for PFDEs with delay. Let us consider the case $0 \le t \le \tau$ in the our model. Since $- \tau  \le t - \tau  \le 0$, then $u(x,t - \tau ) = g(x,t - \tau )$ and therefore Eq (\ref{eq4}) can be rewritten as
\begin{equation}\label{eq10000}
\frac{{\partial u(x,t)}}{{\partial t}} = {K_\alpha }\frac{{{\partial ^\alpha }u(x,t)}}{{\partial {{\left| x \right|}^\alpha }}} + {K_\beta }\frac{{{\partial ^\beta }u(x,t)}}{{\partial {{\left| x \right|}^\beta }}} + \hat g(x,t),\,\,\,\,0 < t \le \tau,
\end{equation}
where $\hat g(x,t) = f(x,t) - {K_\gamma }\dfrac{{{\partial ^\gamma }g(x,t - \tau )}}{{\partial {t^\gamma }}}$. Since $f$ and $\dfrac{{{\partial ^\gamma }g(x,t - \tau )}}{{\partial {t^\gamma }}}$ are continuous functions, $\hat g(x,t)$ is also continuous and so (\ref{eq10000}) is solvable. Thus, we can change initial delay function to control dynamic properties of solutions.\\
We consider the following three fractional advection–dispersion models as the following form
\begin{itemize}
  \item When ${K_\gamma}=0$, Eq (\ref{eq4}) reduces to RFADE \cite{yang2010numerical}.
  \item When ${K_\gamma}=0$ and ${K_\alpha}=0$ , Eq (\ref{eq4}) reduces to Riesz space fractional diffusion equation(RFDE)\cite{yang2010numerical}.
  \item When ${K_\gamma}>0$, Eq (\ref{eq4}) reduces to RFADE with delay (RFADED)\\(our main problem).
\end{itemize}
The fractional advection-dispersion equation is used in groundwater hydrology to model the transport of passive tracers carried by fluid flow in a porous medium \cite{momani2008numerical,liu2003time,huang2005fundamental}. At the second models, physical considerations of a fractional advection–dispersion transport model restrict $0 < \alpha  < 1,\,\,1 < \beta  \le 2$ and we assume ${K_\alpha}>0$ and ${K_\beta \ge 0}$ so that the flow is from left to right. The physical meaning of using homogeneous Dirichlet boundary conditions is that the boundary is set far enough away from an evolving plume such that no significant concentrations reach that boundary \cite{meerschaert2006finite}. In this paper we consider RFADE with delay (${K_\gamma}=1$).\\
Many models with Riesz-space fractional derivatives have recently been studied by many authors.
To approximate the Riesz fractional derivative and use to numerical solving FDEs
and FPDEs, there have existed several analytical and numerical methods. Zhang and Liu \cite{zhang2007fundamental} obtained the analytical solutions of space Riesz and time Caputo fractional partial differential equations. Chen et al. \cite{chen2008fundamental} studied the analytical solution for the space Riesz fractional reaction dispersion equation by using the Laplace and Fourier transform. Yang et al. \cite{yang2010numerical} studied the space Riesz fractional diffusion and advection–dispersion equations and solving these problems on a finite domain by using the L1, L2 methods. Later on, they \cite{zhang2010galerkin} presented two numerical methods for solving the Riesz fractional advectiondispersion equation by using the Galerkin finite element and a backward difference
methods. Shen et al. \cite{shen2014novel} used a second-order numerical scheme for the space Riesz
fractional advection-dispersion equation based on fractional central difference operator. Yang et al. \cite{yang2010numerical} presented three numerical methods for the space Riesz fractional
diffusion and advection-dispersion equations on a finite domain. Celik and Duman
\cite{ccelik2012crank} obtained the numerical solutions of the fractional diffusion equation with the
space Riesz fractional derivative on a finite domain based on the fractional central
difference operator. Recently, Ding et al. \cite{ding2015high} constructed two different fourth-order
numerical schemes for the Riesz derivative and applied to the spatial Riesz fractional diffusion equation. Later on, they \cite{ding2015highII} developed two high-order approximate
schemes (sixth- and eighth-order) for the Riesz fractional derivative and use to solve
the following space Riesz fractional reaction-dispersion equation. Very recently, Ding
and Li \cite{ding2016high} constructed high-order (from 2nd-order to 6th-order) numerical schemes
to approximate the Riesz derivatives and develop three kinds of difference schemes for the Riesz space fractional turbulent diffusion equation. Also, Ding and Li \cite{ding2017high} derived even-order fractional-compact numerical differential formulas for Riesz derivatives and the fourth-order formula applied to solving the Riesz spatial fractional telegraph equation. To study the long-term dynamic behaviors of these models based on Riesz derivatives, analytical solution will play a vital role.
In general, since the fractional integral and fractional derivatives have non-local properties. Therefore, it seems difficult to find the analytical solution of most of these equations. Thus, it is necessary and important to design the numerical schemes to solve these equations.
 Other numerical methods for solving the fractional advection-dispersion equation are investigated in \cite{meerschaert2004finite,liu2004numerical,shen2008fundamental,shen2011numerical,ding2015high,sousa2012second,yang2010numerical}.\\
In this paper, we consider Eqs.(\ref{eq4}-\ref{eq5}). Firstly, we obtain analytical solution of this equation. Secondly, we use fractional backward differential formulas method of second order for $0<\alpha<1$ and shifted Gr\"{u}nwald difference (WSGD) operators for $1 < \beta  \le 2$ to approximate the Riesz space fractional derivative. Next, we obtain the Crank-Nicolson scheme by using the finite difference method for the RFADED. Also, we prove that the Crank-Nicolson scheme is conditionally stable and convergent with the accuracy ${\rm O}({\kappa^2}  + {h^2})$.\\
The paper is organized as follows. In Section 2, we recall some basic definitions of fractional calculus. In Section 3, analytical solution of problem is presented. Fractional backward differential formulas(FBDF) method of second order is introduced in Section 4. In Section 5, shifted Gr\"{u}nwald method is presented. The presented numerical method is proposed and applied to solve the Riesz space fractional advection-dispersion equations with delay (RFADED), stability and convergency of method are proved in Section 6. Finally, some numerical results are given in order to confirm the theoretical analyzes in Section 7.

\section{Preliminaries}
In this section, we will introduce some of the fundamental definitions. Let $\mathbb{C(J,\mathbb{R})}$ denotes the Banach space of all continuous functions from $\mathbb{J} = [0,T]$ into $\mathbb{R}$ with the norm
\begin{equation}\label{eq8}
{\left\| u \right\|_\infty } = \sup \{ \left| {u(t)} \right|:t \in J\} ,\,\,\,\,T > 0.
\end{equation}
${\mathbb{C}^n}(\mathbb{J},\mathbb{R})$ denotes the class of all real valued functions defined on $\mathbb{J} = [0,T]$, $ T>0$ which have continuous $n$th order derivatives.

{\bf Definition 2.1} (\cite{kilbas2006theory}).
The fractional integral of order $\alpha>0$ of the function $f \in {\mathbb{C}}(\mathbb{J},\mathbb{R})$ is defined as
\begin{equation}\label{eq9}
{I^\alpha }f(t) = \frac{1}{{\Gamma (\alpha )}}\int\limits_0^t {\frac{{f(s)}}{{{{(t - s)}^{1 - \alpha }}}}ds,} \,\,\,\,{\mkern 1mu} {\mkern 1mu} {\mkern 1mu} 0 < t < T.
\end{equation}
{\bf Definition 2.2} (\cite{kilbas2006theory}).
The Caputo fractional derivative of order $\alpha>0$ of the function $f \in {\mathbb{C}^n}(\mathbb{J},\mathbb{R})$ is defined as
\begin{equation}\label{eq11*}
{}^C{D^\alpha }f(t) = \left\{ \begin{array}{l}
{I^{n - \alpha }}{D^n}f(t) = \frac{1}{{\Gamma (n - \alpha )}}\int\limits_0^t {\frac{{{f^{(n)}}(s)}}{{{{(t - s)}^{\alpha  - n + 1}}}}ds,} \\
\,\,\,\,\,\,\,\,\,\,\,\,\,\,\,\,\,\,\,\,\,\,\,\,\,\,\,\,\,\,\,\,\,\,\,\,\,\,\,\,\,\,\,\,\,\,\,\,\,\,\,\,\,\,\,\,\,\,n - 1 < \alpha  < n,\,\,n \in N,\\
{f^{(n)}}(t),\,\,\,\,\,\,\,\,\,\,\,\,\,\,\,\,\,\,\,\,\,\,\,\,\,\,\,\,\,\,\,\,\,\,\,\,\,\,\,\,\,\,\,\,\,\,\,\,\,\,\,\,\,\,\,\,\,\,\,\,\,\,\,\,\,\,\,\alpha  = n.
\end{array} \right.
\end{equation}
{\bf Definition 2.3} (\cite{kilbas2006theory,yang2010numerical}).
The Riesz space fractional operator of order $\alpha>0$ of the function $f \in {\mathbb{C}^n}(\mathbb{J},\mathbb{R})$ on a finite domain [0, L] is defined as
\begin{equation}\label{eq6}
\frac{{{\partial ^\zeta }f(x)}}{{\partial {{\left| x \right|}^\zeta }}} =- {( - \Delta )^{\frac{\zeta }{2}}}f(x)= - {c_\zeta }[{}_0^{RL}D_x^\zeta f(x) + {}_x^{RL}D_L^\zeta f(x)],
\end{equation}
where
\begin{equation}\label{eq6}
{c_\zeta } = \frac{1}{{2\cos (\frac{{\pi \zeta }}{2})}},\,\,\,\,0 < \zeta  \le 2,\,\,\zeta  \ne 1
\end{equation}

{\bf Definition 2.4} (\cite{kilbas2006theory}).
 Mittag-leffler functions defined by
\begin{equation}\label{eq1000}
  E_{\alpha ,\beta } (x) = \sum\limits_{k = 0}^\infty  {\frac{{x^k }}
{{\Gamma (\alpha k + \beta )}}} \,\,,\,\,\,x,\beta  \in \mathbb{C},Re(\alpha)>0, \hfill \\
  E_\alpha  (x) = E_{\alpha ,1}.
\end{equation}

{\bf Definition 2.5} (\cite{vcermak2016stability}). The
generalized delay exponential function (of Mittag--Leffler type) is
given by
\begin{equation}\label{eq12}
G_{\alpha ,\beta }^{\lambda ,\tau ,m}(t) = \sum\limits_{j = 0}^\infty  {\left( \begin{array}{l}
	j + m\\
	j
	\end{array} \right)} \frac{{{\lambda ^j}{{(t - (m + j)\tau )}^{\alpha (m + j) + \beta  - 1}}}}{{\Gamma (\alpha (m + j) + \beta )}}H(t - (m + j)\tau ),{\mkern 1mu} {\mkern 1mu} t > 0,
\end{equation}
where $ \lambda  \in \mathbb{C}\,,\,\,\alpha ,\beta ,\tau  \in
\mathbb{R} $ and $ m \in \mathbb{Z}$
 and $H(z)$ is the Heaviside step function.
If $ \lambda  \in \mathbb{C}\,,\,\,\alpha ,\beta ,\tau  \in
\mathbb{R} $ and $ m \in \mathbb{Z}$ then laplace transform of $G_{\alpha ,\beta }^{\lambda ,\tau ,m} (t)$ is:
\begin{equation}\label{eq13}
L(G_{\alpha ,\beta }^{\lambda ,\tau ,m} (t))(s) = \frac{{s^{\alpha
- \beta } \exp \{  - ms\tau \} }} {{(s^\alpha   - \lambda \exp \{
- s\tau \} )^{m + 1} }},\,\,\,\,\,\,\,  s>0.
\end{equation}

\section{Analytical solution of problem}
{\bf Lemma 1.} (\cite{ilic2005numerical})
Suppose the Laplacian ($- \Delta$) has complete set of orthonormal eigenfunctions ${\varphi _n}$ corresponding to eigenvalue $\lambda _n^2$  on a bounded region $\Omega$ i.e. $( - \Delta ){\varphi _n} = \lambda _n^2{\varphi _n}$ on $\Omega$; $B(\Omega) = 0$ on $\partial \Omega$ where $B(\Omega)$ is one of the standard three homogeneous boundary conditions. Let
\begin{equation}\label{eq2000}
{F_\delta } = \left\{ {f = \sum\limits_{n = 1}^\infty  {{f_n}{\varphi _n}} ,\,\,{f_n} =  < f,{\varphi _n} > :\,\,\sum\limits_{n = 1}^\infty  {{{\left| {{f_n}} \right|}^2}\left| \lambda  \right|_n^\delta }  < \infty ,\,\delta  = \max (\alpha ,0)} \right\},
\end{equation}
then for any $f \in {F_\delta },\,{( - \Delta )^{\frac{\alpha }{2}}}:\,\,{F_\delta } \to {L_2}(\Omega )$
is defined by
\begin{equation}\label{eq3000}
{( - \Delta )^{\frac{\alpha }{2}}}f = \sum\limits_{n = 1}^\infty  {{f_n}} {(\lambda _n^2)^{\frac{\alpha }{2}}}{\varphi _n}.
\end{equation}
We consider Eqs. (\ref{eq4}--\ref{eq5}), where ${\mu _1}(t)$ and ${\mu _2}(t)$ are nonzero smooth functions with order--one continuous derivatives. Let
\begin{equation}\label{eq13}
u(x,t) = V(x,t) + W(x,t),
\end{equation}
where
\begin{equation}\label{eq14}
V(x,t) = {\mu _1}(t) +x \frac{{{\mu _2}(t) - {\mu _1}(t)}}{L}.
\end{equation}
By using Eqs. (\ref{eq13}--\ref{eq14}), we can transform the nonhomogeneous condition into a homogeneous boundary condition.
By substituting (\ref{eq13}) into (\ref{eq4}), we have
\begin{equation}\label{eq15}
\begin{array}{l}
\frac{{\partial W(x,t)}}{{\partial t}} +\frac{{{\partial ^\gamma }W(x,t - \tau)}}{{\partial {{t}^\gamma }}}+ {K_\alpha }{( - \Delta )^{\frac{\alpha }{2}}}W(x,t) + {K_\beta }{( - \Delta )^{\frac{\beta }{2}}}W(x,t)= {f_w}(x,t),\,\,t > 0,\\
W(x,t) = {\phi _w}(x,t),\,\,\, - \tau  \le t \le 0,\,\,\,0 \le x \le L,\\
W(0,t) = W(L,t) = 0,\,\,\,\,\,t \ge 0,
\end{array}
\end{equation}
where
\begin{equation}\label{eq16}
\begin{array}{l}
{f_w}(x,t) = f(x,t)-\frac{{{\partial ^\gamma }V(x,t - \tau)}}{{\partial {{t}^\gamma }}} - \frac{{\partial V(x,t)}}{{\partial t}} - {K_\alpha }{( - \Delta )^{\frac{\alpha }{2}}}V(x,t) - {K_\beta }{( - \Delta )^{\frac{\beta }{2}}}V(x,t),\\
{\phi _w}(x,t) = g(x,t) - {\mu _1}(t) - \frac{{{\mu _2}(t) - {\mu _1}(t)}}{L}x.
\end{array}
\end{equation}
We assume that the solution of (\ref{eq15}) has the form:
\begin{equation}\label{eq17}
W(x,t) = X(x)T(t).
\end{equation}
By substituting (\ref{eq17}) into (\ref{eq15}), we obtain the Sturm--Liouville problem with $\lambda >0$ as the following form.
\begin{equation}\label{eq18}
\begin{array}{l}
 - {K_\alpha }{( - \Delta )^{\frac{\alpha }{2}}}X(x) - {K_\beta }{( - \Delta )^{\frac{\beta }{2}}}X(x) + \lambda X(x) = 0,\\
X(0) = 0,\,\,\,\,\,X(L) = 0
\end{array}
 \end{equation}
and the following ODE with delay
\begin{equation}\label{eq19}
\begin{array}{l}
\frac{{dT(t)}}{{dt}} + \frac{{{d^\gamma }T(t - \tau )}}{{d{t^\gamma }}} + \lambda T(t) = 0,\\
T(t) = \rho (x,t),\,\, - \tau  \le t \le 0.
\end{array}
\end{equation}
By using Lemma 1, eigenvalues and corresponding eigenfunctions of the Sturm--Liouville problem (\ref{eq18})--(\ref{eq19}) have the following form
\begin{equation}\label{eq20}
{\lambda _n} = \frac{{{n^2}{\pi ^2}}}{{{L^2}}},\,\,\,\,{X_n}(x) = \sin (\frac{{n\pi }}{L}x),\,\,\,n = 1,2,...\cdot
\end{equation}
Therefore, we set
\begin{equation}\label{eq21}
W(x,t) = \sum\limits_{n = 1}^\infty  {{B_n}(t)\sin (} \frac{{n\pi }}{L}x).
\end{equation}
By substituting (\ref{eq21}) into (\ref{eq15}), we can write
\begin{equation}\label{eq22}
\begin{array}{*{20}{l}}
{\frac{{d{B_n}(t)}}{{dt}} + \frac{{{d^\gamma }{B_n}(t - \tau )}}{{d{t^\gamma }}} + {\Theta _n}{B_n}(t) +  = {f_{wn}}(t),}\\
{{B_n}(t) = \varphi (t),}
\end{array}
\end{equation}
where ${\Theta _n} = {K_\alpha }{\lambda _n}^{\frac{\alpha }{2}} + {K_\beta }{\lambda _n}^{\frac{\beta }{2}}$ and
\begin{equation}\label{eq23}
\begin{array}{l}
{f_{wn}}(t) = \frac{2}{L}\int_0^L {{f_w}(x,t)\sin (\frac{{n\pi }}{L}x)dx} ,\\
{f_w}(x,t) = \sum\limits_{n = 1}^\infty  {{f_{wn}}(t)\sin (} \frac{{n\pi }}{L}x),\\
\varphi (t) = \frac{2}{L}\int_0^L {{\phi _w}(x,t)\sin (\frac{{n\pi }}{L}x)dx}.\\
\end{array}
\end{equation}
Firstly, we have
\begin{equation}\label{eq23*}
\begin{array}{l}
L\{ \frac{{{d^\gamma }{B_n}(t - \tau )}}{{d{t^\gamma }}},s\}  = {s^\gamma }L\{ {B_n}(t - \tau ),s\}  - {s^{\gamma  - 1}}{B_n}( - \tau )\\
\,\,\,\,\,\,\,\,\,\,\,\,\,\,\,\,\,\,\,\,\,\,\,\,\,\,\,\,\,\,\,\,\,\,\,\,\,\,\,\,\, = \,{s^\gamma }\{ {e^{ - s\tau }}\overline {{B_n}(s)}  + {e^{ - s\tau }}\int\limits_{ - \tau }^0 {{e^{ - s\tau }}{B_n}(\nu )d\nu } \}  - {s^{\gamma  - 1}}\varphi ( - \tau ),
\end{array}
\end{equation}
where $L$ is the Laplace operator. By using Laplace transform for (\ref{eq22}) and Eq. \ref{eq23*}, we can write
\begin{equation}\label{eq24}
\begin{array}{l}
\overline {{B_n}(s)}  = \frac{{{s^{ - \gamma }}{{\bar f}_{wn}}(s)}}{{{s^{1 - \gamma }} + {e^{ - \tau s}} + {\Theta _n}{s^{ - \gamma }}}} + \frac{{{s^{ - \gamma }}\varphi (0)}}{{{s^{1 - \gamma }} + {e^{ - \tau s}} + {\Theta _n}{s^{ - \gamma }}}} + \frac{{{s^{ - 1}}\varphi ( - \tau )}}{{{s^{1 - \gamma }} + {e^{ - \tau s}} + {\Theta _n}{s^{ - \gamma }}}}\\
\,\,\,\,\,\,\,\,\,\,\,\,\, - \frac{{{e^{ - \tau s}}\int_{ - \tau }^0 {{e^{ - s\nu }}} \varphi (\nu )d\nu }}{{{s^{1 - \gamma }} + {e^{ - \tau s}} + {\Theta _n}{s^{ - \gamma }}}},
\end{array}
\end{equation}
where
\begin{equation}\label{eq25}
\overline {{B_n}(s)}  = L({B_n}(t)),{\mkern 1mu} {\mkern 1mu} {{\bar f}_{wn}}(s) = L({f_{wn}}(t)).
\end{equation}
Also, we can write
\begin{equation}\label{eq26}
\begin{array}{l}
\frac{1}{{{s^{1 - \gamma }} + {e^{ - \tau s}} + {\Theta _n}{s^{ - \gamma }}}} = \frac{1}{{{\Theta _n}{s^{ - \gamma }}}}\frac{{{\Theta _n}{s^{ - \gamma }}}}{{{s^{1 - \gamma }} + {e^{ - \tau s}}}}\frac{1}{{1 + \frac{{{\Theta _n}{s^{ - \gamma }}}}{{{s^{1 - \gamma }} + {e^{ - \tau s}}}}}}{\mkern 1mu} {\mkern 1mu} {\mkern 1mu} {\mkern 1mu} {\mkern 1mu} {\mkern 1mu} {\mkern 1mu} \\
\,\,\,\,\,\,\,\,\,\,\,\,\,\,\,\,\,\,\,\,\,\,\,\,\,\,\,\,\,\,\,\,\,\,\,\,\,\,\, = \sum\limits_{k = 0}^\infty  {{{( - {\Theta _n})}^k}{e^{k\tau s}}\frac{{{s^{1 - \gamma  - (k\gamma  - \gamma  + 1)}}{e^{ - k\tau s}}}}{{{{({s^{1 - \gamma }} + {e^{ - \tau s}})}^{k + 1}}}}} .
\end{array}
\end{equation}
Therefore
\begin{equation}\label{eq27}
{L^{ - 1}}(\frac{1}{{{s^{1 - \gamma }} + {e^{ - \tau s}} + {\Theta _n}{s^{ - \gamma }}}}) = \sum\limits_{k = 0}^\infty  {{{( - {\Theta _n})}^k}G_{1 - \gamma ,(k - 1)\gamma  + 1}^{ - 1,\tau ,k}(t + k\tau )} H(t + k\tau ).
\end{equation}
We assume that
\begin{equation}\label{eq29}
Z(s) = \int_{ - \tau }^0 {{e^{ - s\nu }}} \varphi (\nu )d\nu
\end{equation}
and
\begin{equation}\label{eq30}
{L^{ - 1}}(Z(s)) = z(t).
\end{equation}
Therefore, we have
\begin{equation}\label{eq28}
\begin{array}{l}
{L^{ - 1}}(\frac{{{s^{ - \gamma }}\varphi (0)}}{{{s^{1 - \gamma }} + {e^{ - \tau s}} + {\Theta _n}{s^{ - \gamma }}}}) = \varphi (0)\sum\limits_{k = 0}^\infty  {{{( - {\Theta _n})}^k}G_{1 - \gamma ,k\gamma  + 1}^{ - 1,\tau ,k}(t + k\tau )} H(t + k\tau ),\\
{L^{ - 1}}(\frac{{{s^{ - 1}}\varphi ( - \tau )}}{{{s^{1 - \gamma }} + {e^{ - \tau s}} + {\Theta _n}{s^{ - \gamma }}}}) = \varphi ( - \tau )\sum\limits_{k = 0}^\infty  {{{( - {\Theta _n})}^k}G_{1 - \gamma ,(k - 1)\gamma  + 2}^{ - 1,\tau ,k}(t + k\tau )} H(t + k\tau ),\\
{L^{ - 1}}(\frac{{{s^{ - \gamma }}{{\bar f}_{wn}}(s)}}{{{s^{1 - \gamma }} + {e^{ - \tau s}} + {\Theta _n}{s^{ - \gamma }}}})\\
\,\,\,\,\,\,\,\,\,\,\,\,\,\,\,\,\,\,\,\,\,\,\,\,\,\,\,\,\,\,\,\,\,\,\,\,\,\,\,\,\,\,\,\,\, = \sum\limits_{k = 0}^\infty  {{{( - {\Theta _n})}^k}\int_0^t {G_{1 - \gamma ,k\gamma  + 1}^{ - 1,\tau ,k}(t - p + k\tau )H(t - p + k\tau ){f_{wn}}(p)dp} }
\end{array}
\end{equation}
and
\begin{equation}\label{eq28*}
\begin{array}{l}
{L^{ - 1}}(\frac{{{e^{ - \tau s}}\int_{ - \tau }^0 {{e^{ - s\nu }}} \varphi (\nu )d\nu }}{{{s^{1 - \gamma }} + {e^{ - \tau s}} + {\Theta _n}{s^{ - \gamma }}}})\\
\,\,\,\,\,\,\,\,\,\,\,\,\,\,\,\, = \sum\limits_{k = 0}^\infty  {{{( - {\Theta _n})}^k}\int_0^t {G_{1 - \gamma ,(k - 1)\gamma  + 1}^{ - 1,\tau ,k}(t - p + (k - 1)\tau )H(t - p + (k - 1)\tau )z(p)dp} } .
\end{array}
\end{equation}
Therefore, we have
\begin{equation}\label{eq32}
\begin{array}{l}
{B_n}(t) = \sum\limits_{k = 0}^\infty  {{{( - {\Theta _n})}^k}[\int_0^t {(G_{1 - \gamma ,k\gamma  + 1}^{ - 1,\tau ,k}({\rm{t  -  p}} + k\tau )H({\rm{t  -  p}} + k\tau )f(p)} } \\
\,\,\,\,\,\,\,\,\,\,\,\,\,\,\, - G_{1 - \gamma ,(k - 1)\gamma  + 1}^{ - 1,\tau ,k}({\rm{t  -  p}} + (k - 1)\tau )H({\rm{t  -  p}} + (k - 1)\tau )z(p))dp\\
\,\,\,\,\,\,\,\,\,\,\,\,\,\, + \,\varphi (0)G_{1 - \gamma ,k\gamma  + 1}^{ - 1,\tau ,k}(t + k\tau )H(t + k\tau )\\
\,\,\,\,\,\,\,\,\,\,\,\,\,\, + \varphi ( - \tau )G_{1 - \gamma ,(k - 1)\gamma  + 2}^{ - 1,\tau ,k}(t + k\tau )H(t + k\tau )].
\end{array}
\end{equation}
Finally, the analytical solution of (\ref{eq15})--(\ref{eq16}), is
\begin{equation}\label{eq33}
\begin{array}{l}
u(x,t) = {\mu _1}(t) + \frac{{{\mu _2}(t) - {\mu _1}(t)}}{L}x{\mkern 1mu} \\
\,\,\,\,\,\,\,\,\,\,\,\,\,\,\,\, + \sum\limits_{n = 1}^\infty  {\sum\limits_{k = 0}^\infty  {{{( - {\Theta _n})}^k}[\int_0^t {(G_{1 - \gamma ,k\gamma  + 1}^{ - 1,\tau ,k}({\rm{t  -  p}} + k\tau )H({\rm{t  -  p}} + k\tau )f(p)} } } \\
\,\,\,\,\,\,\,\,\,\,\,\,\,\,\,\,\, - G_{1 - \gamma ,(k - 1)\gamma  + 1}^{ - 1,\tau ,k}({\rm{t  -  p}} + (k - 1)\tau )H({\rm{t  -  p}} + (k - 1)\tau )z(p))dp\\
\,\,\,\,\,\,\,\,\,\,\,\,\,\,\,\,\, + \,\varphi (0)G_{1 - \gamma ,k\gamma  + 1}^{ - 1,\tau ,k}(t + k\tau )H(t + k\tau )\\
\,\,\,\,\,\,\,\,\,\,\,\,\,\,\,\,\, + \varphi ( - \tau )G_{1 - \gamma ,(k - 1)\gamma  + 2}^{ - 1,\tau ,k}(t + k\tau )H(t + k\tau )]sin(\frac{{n\pi }}{L}x).
\end{array}
\end{equation}

\section{Fractional Backward Differential Formulas(FBDF) method of second order}
\subsection{Discrete convolution quadratures}
 We consider the initial value problem
 \begin{equation}\label{eq6*}
 {}_{t_0 }^C D_t^\alpha  y(t) = f(t),\,\,\,\,\,\,0\, < \alpha  < 1,
 \hfill
 \end{equation}
 \[
 y (t_0 ) = y_{0},
 \]
 where $f$ is a sufficiently smooth function. Problem (\ref{eq6*})
 can be written in the form (see \cite{kilbas2006theory})
 \begin{equation}\label{eq7*}
 y(t) = y_0  + \frac{1} {{\Gamma (\alpha )}}\int_{t_0 }^t {(t -
 	\tau )^{\alpha  - 1} f(\tau)d\tau }.
 \end{equation}
 Integral (\ref{eq7*}), may be approximated by

 \begin{equation}\label{eq88}
 I_h^\alpha  f(t_n ) = h^\alpha  \sum\limits_{j = 0}^n {\omega _j
 	f(t_{n - j} ) + h^\alpha  \sum\limits_{j = 0}^s {\omega _{n,j}
 		f(t_j )} } ,\,\,t_j  = t_0  + jh,
 \end{equation}

 by means of Eqs. (\ref{eq7*}) and (\ref{eq88}), we obtain
 \begin{equation}\label{eq9*}
 y_n  = y_0  + h^\alpha  \sum\limits_{j = 0}^n {\omega _j f_{n - j}
 	+ h^\alpha  \sum\limits_{j = 0}^s {\omega _{n,j} f_j } },
 \end{equation}
 where $\omega _j$ are coefficients in the formal power series $ \omega (z) = \sum\limits_{j = 0}^\infty  {\omega _j z^j } $ and  $\omega _{n,j}$ are starting quadrature weights (\cite{lubich1986discretized,henrici1962discrete}).

 \subsection{Fractional linear multistep methods}
 Fractional linear multistep methods (FLMMs)
 for the approximate solution of Eq. (\ref{eq6*})
 is as the following form (see\cite{lubich1986discretized})
 \begin{equation}\label{eq008}
 \sum\limits_{j = 0}^n {\alpha _j y_{n - j}  = h^{\alpha} \sum\limits_{j =
 		0}^n {\beta _j f_{n - j} } },
 \end{equation}
 which
 \[
 \rho (\xi ) = \sum\limits_{j = 0}^n {\alpha _j \xi ^j }
 ,\,\,\,\,\,\,\sigma (\xi ) = \sum\limits_{j = 0}^n {\beta _j \xi
 	^j },
 \]
 where $ \rho$ and $\sigma$ denote the generating polynomials of the method (see \cite{henrici1962discrete}).\\
 As shown in \cite{w1979linear},  FLMMs
 (\cite{lubich1986discretized,henrici1962discrete}) can be
 equivalently reformulated as Eq. (\ref{eq9*}) with weights $\omega_j$ obtained as the coefficients in the FPS of the fractional-order power of the generating function (see \cite{lubich1986discretized}):
 \begin{equation}\label{eq001}
 \omega ^{(\alpha )} (\xi) = \left( \frac{{\sigma
 		(\frac{1} {\xi})}} {{\rho (\frac{1} {\xi})}}\right) ^\alpha= \sum\limits_{j = 0}^\infty  {\omega _j
 	\xi^j }.
 \end{equation}
 Methods of this kind, named as FLMMs, when applied to (\ref{eq6*}) we will have Eq. (\ref{eq9*})
 where weights are given in (\ref{eq001}).

 \subsection{The second order FBDF}
 The n-step backward differentiation formulas is the implicit linear multistep formula with ${\beta _0} = {\beta _1} = ... = {\beta _{n - 1}} =
 0$ in (\ref{eq008}). The second order BDF formulas for ODEs is given by
 \[
3y_{n + 1}  - 4y_n  + y_{n - 1}  = 2hf_{n + 1}.
\]
 We take
 $
\rho (\xi ) = 3\xi^2  - 4\xi + 1$ and $ \sigma (\xi ) =
2\xi^2.$ \\
 By using Eq. (\ref{eq001}), we have
\[
\omega ^{(\alpha )} (\xi) = (\frac{3} {2} - 2\xi + \frac{1}
{2}\xi^2 )^{\alpha }  = \sum\limits_{j = 0}^\infty  {\omega _j
\xi^j }.
\]
 By means of Miller's recurrence (\cite{zeilberger1995jcp}), we can obtain
 \begin{equation}\label{eq12*}
	{\omega _0} = 1\,,\,\,\,{\omega _1} = -\frac{4}{3}\alpha ,\,\,\\
	\omega_k = \frac{4} {3}(1 - \frac{{\alpha  + 1}} {k})\omega _{k -
		1}  + \frac{1} {3}(\frac{{2(1 + \alpha )}} {k} - 1)\omega _{k - 2}
	,\,\,\,k = 2,3,\cdots,
\end{equation}
 the weights of fractional backward differentiation formulas of second order (FBDF2) will be as
 \[ \varpi _k  = (\frac{3} {2})^\alpha\omega _k ,\,k = 0,1,\cdots.\]

 \subsection{Numerical method (FBDF2)}
 By using (\ref{eq40000}) and
 (\ref{eq11*}) in (\ref{eq6*}) at $ t = t_j $,  we have
 \begin{equation}
 	{}_{t_0 }^C D_t^\alpha  y(t)= {}_{t_0 }^{RL} D_t^\alpha  (y(t) - y(t_0 )).
 \end{equation}
 Then we can write
 \begin{equation}
  {}_{t_0 }^{RL} D_t^\alpha  y(t_j ) - \frac{{t_j ^{ - \alpha }
 	}} {{\Gamma (1 - \alpha )}}y(t_0 ) = f(t_j),\,\,\,0 <
 	\alpha < 1.
\end{equation}
 We know Gr\"{u}nwald approximation is as the following form (\cite{podlubny1999fractional})
 \begin{equation}\label{eq013}
 {}_{t_0 }^{RL} D_t^\alpha  y(t_j ) \approx h^{ - \alpha }
 \sum\limits_{k = 0}^j {\varpi _k y(} t_j  - kh),
 \end{equation}
 where
 \begin{equation}\label{eq013*}
 (1 - \xi)^\alpha   = \sum\limits_{k = 0}^\infty  {\varpi _k \xi^k
 }.
 \end{equation}
 By using (\ref{eq013}), we can rewrite
 \begin{eqnarray*}
 	{}_{t_0 }^C D_t^\alpha  y(t) \approx h^{ - \alpha } \sum\limits_{k
 		= 0}^j {\varpi _k y(} t_j  - kh) - \frac{{t_j ^{ - \alpha } }}
 	{{\Gamma (1 - \alpha )}}y(t_0 ) = f(t_j),\,\,\,0 < \alpha
 	< 1.
 \end{eqnarray*}
 If we take
 \begin{equation}\label{eq013**}
 \frac{{j^{ - \alpha } }} {{\Gamma (1 - \alpha )}} = b_j ,\,y_{j -
 	k}  = y(t_j  - kh),\,f_j  = f(t_j),
 \end{equation}
 we obtain
 \begin{equation}\label{eq014}
 \sum\limits_{k = 0}^j {\varpi _k y_{j - k} }  - b_j y_0  =
 h^\alpha  f_j.
 \end{equation}
 Therefore FBDF2 is defined as
 \begin{equation}\label{eq015}
 \sum\limits_{k = 0}^j {\varpi _k y_{j - k} }  - b_j y_0  =
 h^\alpha  f_j,
 \end{equation}
 with weights
 \begin{equation}\label{eq016}
\begin{array}{l}
{\omega _0} = 1{\mkern 1mu} ,{\mkern 1mu} {\mkern 1mu} \\
{\mkern 1mu} {\omega _1} =  - \frac{4}{3}\alpha ,{\mkern 1mu} {\mkern 1mu} \\
{\omega _k} = \frac{4}{3}(1 - \frac{{\alpha  + 1}}{k}){\omega _{k - 1}} + \frac{1}{3}(\frac{{2(1 + \alpha )}}{k} - 1){\omega _{k - 2}},\\
{\varpi _k} = {(\frac{3}{2})^\alpha }{\omega _k},\,\,\,k = 0,1, \cdots .
\end{array}
\end{equation}

\subsection{Analysis of linear stability and consistency}
In order to study consistency and linear stability of FBDF
consider the linear test problem
\begin{equation}\label{eq44}
{}_{t_0 }^C D_t^\alpha  y(t) = \lambda y(t),\,\,\,\lambda \in
C,\,\,0 < \alpha  < 1.
\end{equation}
{\bf Theorem 1.} (\cite{garrappa2013trapezoidal} ) The
steady-state $ y=0 $ of (\ref{eq44}) is stable if and only if
\begin{equation}\label{eq45}
\lambda  \in \left\{ {\upsilon  \in C:\,\left| {\arg (\upsilon )}
	\right| > \alpha \frac{\pi } {2}} \right\}.
\end{equation}

\subsubsection{Consistency of FBDF2}
{\bf Theorem 2.} The FBDF2 method is consistent of second order.\\
{\bf Proof.} A convolution quadrature $ \omega^{(\alpha)} $ is consistent of order $p$ if
	\[{h^{ - \alpha }}{\omega ^{(\alpha )}}({e^{ - h}}) = 1 + O({h^p}),\](see \cite{lubich1986discretized}). Thus for FBDF method of second order, we have
\[
\omega^{(\alpha)}(z) = \left( \frac{3} {2} - 2z + \frac{1} {2}z^2 \right) ^\alpha.
\]
Therefore we can write
\[
h^{ - \alpha } \omega^{(\alpha)} (e^{ - h} )= \left( \frac{{h - \frac{1}
		{3}h^3  + ...}} {h}\right) ^\alpha   = \left( 1 - \frac{1} {3}h^2 \right) ^\alpha = 1
+ O(h^2 ).
\]
Thus, the FBDF2 method is consistent of second order.
\subsubsection{Stability of FBDF2}
{\bf Lemma 2.} For weights of FBDF2 we have $ 1 - \varepsilon  <
	\frac{{\omega _n }} {{\omega _{n - 1} }} < 1$  for $n\geq4, 0 <
	\varepsilon  <\frac{1} {3}$.

{\bf Proof.} The proof is by induction on $n$. The relation is true for $n=4$. Suppose $ 1 - \varepsilon  <
	\frac{{\omega _{n - 1} }} {{\omega _{n - 2} }} < 1\,,\,\,\,\,\,0 <
	\varepsilon  < \frac{1} {3}$ is true. We will prove it for $n$. If $ \varsigma _n  = \frac{{\alpha
			+ 1}} {n}$ and
	$n\geq4$ then we have\\
	\[
	1 - \varsigma _n  > 0\,,\,\,1 - 2\varsigma _n  > 0
	\]
	and
	\begin{equation}\label{eq22*}
		1 < \frac{{\omega _{n - 2} }} {{\omega _{n - 1} }} < \frac{1} {{1
				- \varepsilon }}\,\,\, \to \,\, - \frac{1} {{1 - \varepsilon }} <
		- \frac{{\omega _{n - 2} }} {{\omega _{n - 1} }} <  - 1.
	\end{equation}\\
	By using (\ref{eq22*}), we can rewrite
	\[
	\frac{4} {3}(1 -
	\varsigma _n ) - \frac{1} {{1 - \varepsilon }}\frac{1} {3}(1 -
	2\varsigma _n ) < \frac{4} {3}(1 - \varsigma _n ) - \frac{1} {3}(1
	- 2\varsigma _n )\frac{{\omega _{n - 2} }} {{\omega _{n - 1} }} <
	\frac{4} {3}(1 - \varsigma _n ) - \frac{1} {3}(1 - 2\varsigma _n
	).
	\]
	Since $ 0 < \varsigma _n  < \frac{1} {2}$ we have $ \frac{4} {3}(1
	- \varsigma _n ) - \frac{1} {{1 - \varepsilon }}\frac{1} {3}(1 -
	2\varsigma _n ) > 1 - \frac{{1 - 2\varepsilon }} {{3 -
			3\varepsilon }}$. Now put $ \varepsilon ^*  = \frac{{1 - 2\varepsilon
	}} {{3 - 3\varepsilon }}\,$ then from $ 0 < \varepsilon ^*  < \frac{1}
	{3}$ we have
	\[
	\frac{4} {3}(1 - \varsigma _n ) - \frac{1} {{1 - \varepsilon
	}}\frac{1} {3}(1 - 2\varsigma _n ) > 1 - \varepsilon ^*.
	\]
	Finally
	\[
	1 - \varepsilon ^*  < \frac{4} {3}(1 - \varsigma _n ) - \frac{1}
	{3}(1 - 2\varsigma _n )\frac{{\omega _{n - 2} }} {{\omega _{n - 1}
	}} < 1,
	\]
	from the above relation we obtain
	\[
	1 - \varepsilon  < \frac{{\omega _n }} {{\omega _{n - 1} }} <
	1\,,\,\,\,\,\,0 < \varepsilon  <  \frac{1} {3}\,\,,\,n =
	4,5,\cdots.
	\]

{\bf Lemma 3.}
	For weights of FBDF2 and $0<\alpha<1$ we have
	\begin{equation}\label{eq47}
		\omega _0  > 0\,,\,\,\omega _j  < 0\,,\,j = 4,5\cdots.
	\end{equation}
	\[
	 \left|
	{\omega _j } \right| < \left| {\omega _{j - 1} } \right| < \omega
	_0 \,,\,j = 4,5,\cdots
	\]
	
	\[
	\sum\limits_{k = 0}^\infty  {\omega _k  = 0,\,\,\,\,\,\,}
	\,\sum\limits_{k = 0}^{m} {\omega _k } > 0,\,\,\,\,m >3.
	\]
{\bf Proof.} Proof is by induction on $n$ by using Lemma 1. The relation is true for $j=4$. It is obvious that
	\[
	\omega _0  = (\frac{3} {2})^\alpha   > 0\,,\,\,\omega _4  =  -
	\frac{1} {{486}}(\frac{3} {2})^\alpha  \alpha (1 - \alpha )(64(1 -
	\alpha )^2  + 48(1 - \alpha ) + 11) < 0,
	\]
	now suppose that $ \omega _j  < 0\,,\,j = 5,6,...n-1 $. We can
	write
	\[
	\omega _n  = \omega _{n - 1} (\frac{4} {3}(1 - \varsigma _n ) -
	\frac{1} {3}(1 - 2\varsigma _n )\frac{{\omega _{n - 2} }} {{\omega
			_{n - 1} }}).
	\]
	Suppose $ \vartheta  = \frac{4} {3}(1 - \varsigma _n ) - \frac{1}
	{3}(1 - 2\varsigma _n )\frac{{\omega _{n - 2} }} {{\omega _{n - 1}
	}}. $  Since $ \omega _{n - 1}  < 0$ and by using Lemma 1, we have
	\[\begin{array}{l}
	\frac{4}{3}(1 - {\varsigma _n}) - \frac{1}{{1 - \varepsilon }}\frac{1}{3}(1 - 2{\varsigma _n}) < \vartheta  < \frac{4}{3}(1 - {\varsigma _n}) - \frac{1}{3}(1 - 2{\varsigma _n}),\\
	\frac{4}{3}(1 - {\varsigma _n}) - \frac{1}{{1 - \varepsilon
	}}\frac{1}{3}(1 - 2{\varsigma _n}) > 1 - {\varepsilon
		^*},\,\,\,\,\,0 < {\varepsilon ^*} < \frac{1}{3},
	\end{array}\]
	therefore
	\[\frac{4}{3}(1 - {\varsigma _n}) - \frac{1}{{1 - \varepsilon }}\frac{1}{3}(1 - 2{\varsigma _n}) > 0.\]
	Then $\vartheta  > 0$. Thus, we have
	\[
	\omega _n  = \vartheta \omega _{n - 1}  < 0.
	\]
	By means of  Lemma 1 we can write
	\[
	\omega _n  < 0\,\, \to \omega _{n - 1}  < \omega _n \, \to \left|
	{\omega _n } \right| < \left| {\omega _{n - 1} } \right| < \omega
	_0 \,,\,n = 4,5,\cdots.
	\]
	Therefore, we have
	\[
	\sum\limits_{k = 0}^\infty  {\omega _k  = 0\, \to \,}
	\,\sum\limits_{k = 0}^{j - 1} {\omega _k }  =  - (\omega _j  +
	\omega _{j + 1}  + ....) > 0\, \to b_j  > 0.
	\]\\

\begin{figure}[h]
	\begin{center}
		\includegraphics[width=6.8cm]{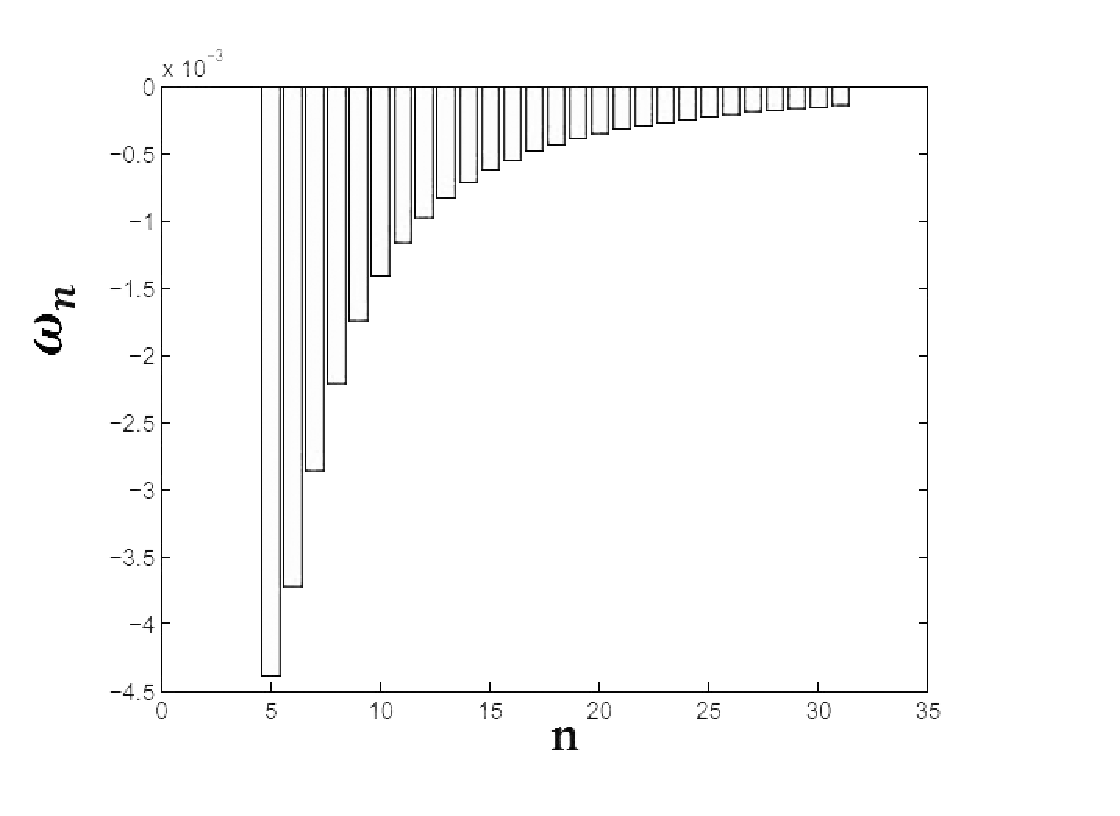}\includegraphics[width=6.3cm]{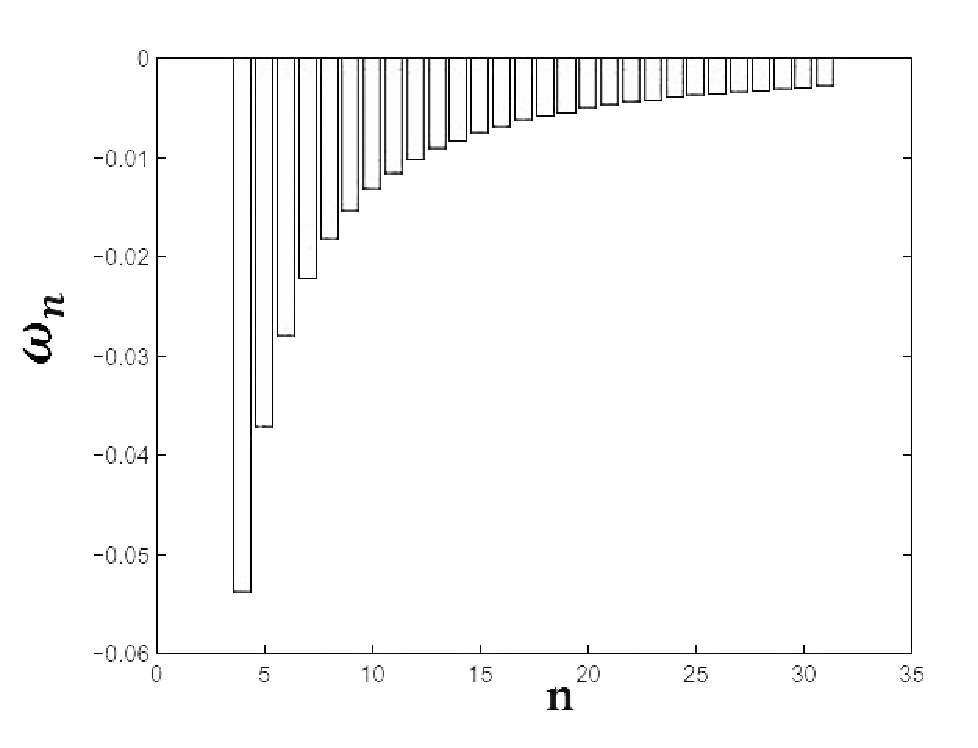}
		\caption{{\small Weights of $FBDF2$ for $\alpha=0.3(left)$ and $\alpha=0.7(right)$.}
		}\label{1}
	\end{center}
\end{figure}

{\bf Theorem 3.} Numerical solution of problem (\ref{eq44}) obtained by FBDF2 is unconditionally stable.\\
{\bf Proof.} By using FBDF2 method for (\ref{eq44}), the matrix form of the problem will be as the following form
\begin{equation}\label{eq48}
\left( {\begin{array}{*{20}{c}}
{{\varpi _0} - \lambda {h^\alpha }}&0& \cdots &{\begin{array}{*{20}{c}}
{}&0
\end{array}}\\
{{\varpi _1}}&{{\varpi _0} - \lambda {h^\alpha }}&0&{\begin{array}{*{20}{c}}
 \cdots &0
\end{array}}\\
 \vdots & \vdots & \vdots &{\begin{array}{*{20}{c}}
{}& \vdots
\end{array}}\\
{{\varpi _{n - 1}}}&{{\varpi _{n - 2}}}& \cdots &{{\varpi _0} - \lambda {h^\alpha }}
\end{array}} \right)\left( \begin{array}{l}
{y_1}\\
{y_2}\\
 \vdots \\
{y_n}
\end{array} \right) = {y_0}\left( \begin{array}{l}
{b_1}\\
{b_2}\\
 \vdots \\
{b_n}
\end{array} \right)
\end{equation}
If $\lambda<0$ then by means of Lemma's 2-3, coefficient matrix is strictly diagonally dominant. Therefore numerical solution of problem obtained by the FBDF2 method will be unconditionally stable.
\subsubsection{Convergent of FBDF2}

{\bf Theorem 4.} (\cite{lubich1986discretized}) A convolution quadrature $ \omega^{(\alpha)} $ is convergent of order $p$ if and only if it is stable and consistent of order $p$.\\
{\bf Theorem 5.} The FBDF2 method is convergent of second order.\\
{\bf Proof.} Since FBDF2 is consistence of second order and stable by means of Theorem 4, therefore the present method is convergent of second order.

\section{Shifted Gr\"{u}nwald method}
The Shifted Gr\"{u}nwald difference operators is as  the following form
\begin{equation}\label{eq49}
\begin{array}{l}
M_{h,p}^\zeta u(x) = {h^{ - \zeta }}\sum\limits_{k = 0}^\infty  {\omega _k^{(\zeta )}u(x - (k - p)h)} ,\\
N_{h,q}^\zeta u(x) = {h^{ - \gamma }}\sum\limits_{k = 0}^\infty  {\omega _k^{(\zeta )}u(x + (k - q)h)}.
\end{array}
\end{equation}
Therefore, we have
\begin{equation}\label{eq50}
\begin{array}{l}
M_{h,p}^\zeta u(x) = {}_{ - \infty }^{}D_x^\zeta u(x) + {\rm O}(h),\\
N_{h,p}^\zeta u(x) = {}_x^{}D_{ + \infty }^\zeta u(x) + {\rm O}(h),
\end{array}
\end{equation}
where $p,q \in Z$ and $\omega _k^{(\zeta )} = {( - 1)^k}\left( \begin{array}{l}
\zeta \\
k
\end{array} \right)$.\\
{\bf Lemma 4.} (\cite{tian2015class}) Suppose that $1 < \zeta  \le 2$, then the coefficients $\omega _k^{(\zeta )}$  satisfy
\begin{equation}\label{eq51}
\begin{array}{l}
\omega _0^{(\zeta )} = 1,\,\,\,\omega _1^{(\zeta )} =  - \zeta ,\,\,\omega _2^{(\zeta )} = \frac{{\zeta (\zeta  - 1)}}{2},\\
1 \ge \omega _2^{(\zeta )} \ge \omega _3^{(\zeta )} \ge ... \ge 0,\\
\sum\limits_{k = 0}^\infty  {\omega _k^{(\zeta )} = 0,\,\,\,\,\sum\limits_{k = 0}^m {\omega _k^{(\zeta )}} }  < 0,\,\,\,m \ge 1.
\end{array}
\end{equation}
{\bf Theorem 6.} (\cite{tian2015class}) Let $1 < \zeta  \le 2$ and $u \in {L^1}(R),\,\,{}_{ - \infty }^{}D_x^\zeta u,\,\,{}_x^{}D_{ + \infty }^\zeta u$ and their Fourier transforms belong to ${L^1}(R)$ and define the weighted and shifted Gr\"{u}nwald difference (WSGD) operator by
\begin{equation}\label{eq52}
\begin{array}{l}
{}_L^{}D_{h,p,q}^\zeta u(x) = \frac{{\zeta  - 2q}}{{2(p - q)}}M_{h,p}^\zeta u(x) + \frac{{2p - \zeta }}{{2(p - q)}}M_{h,q}^\zeta u(x),\\
{}_R^{}D_{h,p,q}^\zeta u(x) = \frac{{\zeta  - 2q}}{{2(p - q)}}N_{h,p}^\zeta u(x) + \frac{{2p - \zeta }}{{2(p - q)}}N_{h,q}^\zeta u(x),
\end{array}
\end{equation}
then we have
\begin{equation}\label{eq53}
\begin{array}{l}
{}_L^{}D_{h,p,q}^\zeta u(x) = {}_{ - \infty }^{}D_x^\zeta u(x) + {\rm O}({h^2}),\\
{}_R^{}D_{h,p,q}^\zeta u(x) = {}_x^{}D_{ + \infty }^\zeta u(x) + {\rm O}({h^2}),
\end{array}
\end{equation}
where $p$ and $q$ are integers and $p \ne q$.\\
If  for well defined function $u(x)$ on the bounded interval $[a, b]$, we have $u(a) = 0$ or $u(b) = 0$, the function $u(x)$ can be zero extended for $x < a$ or $x > b$. Therefore the left and right Riemann--Liouville fractional derivatives of $u(x)$ at each point $x$ can be approximated by the WSGD operators with second order accuracy
\begin{equation}\label{eq54}
\begin{array}{l}
{}_a^{}D_x^\zeta u(x) = {\eta _1}{h^{ - \zeta }}\sum\limits_{k = 0}^{\left[ {\frac{{x - a}}{h}} \right] + p} {\omega _k^{(\zeta )}u(x - (k - p)} h)\\
\,\,\,\,\,\,\,\,\,\,\,\,\,\,\,\,\,\,\,\,\,\,\,\,\,\,\,\,\,\,\,\,\,\,\,\,\,\,\,\,\,\,\,\, + {\eta _2}{h^{ - \zeta }}\sum\limits_{k = 0}^{\left[ {\frac{{x - a}}{h}} \right] + q} {\omega _k^{(\zeta )}u(x - (k - q)} h) + {\rm O}({h^2}),\\
{}_x^{}D_b^\zeta u(x) = {\eta _1}{h^{ - \zeta }}\sum\limits_{k = 0}^{\left[ {\frac{{b - x}}{h}} \right] + p} {\omega _k^{(\zeta )}u(x + (k - p)} h)\\
\,\,\,\,\,\,\,\,\,\,\,\,\,\,\,\,\,\,\,\,\,\,\,\,\,\,\,\,\,\,\,\,\,\,\,\,\,\,\,\,\,\,\,\, + {\eta _2}{h^{ - \zeta }}\sum\limits_{k = 0}^{\left[ {\frac{{b - x}}{h}} \right] + q} {\omega _k^{(\zeta )}u(x + (k - q)} h) + {\rm O}({h^2}),
\end{array}
\end{equation}
where ${\eta _1} = \frac{{\zeta  - 2q}}{{2(p - q)}},\,{\eta _2} = \,\frac{{2p - \zeta }}{{2(p - q)}}$.\\
For $1 < \zeta \le 2$, if we take $(p, q) = (1, 0)$, therefore Eq (\ref{eq54}) on the domain $[0, L]$, will be as the following form
\begin{equation}\label{eq55}
\begin{array}{l}
{}_0^{}D_x^\zeta u({x_i}) = {h^{ - \zeta }}\sum\limits_{k = 0}^{i + 1} {\vartheta _k^{(\zeta )}u({x_{i - k + 1}}} ) + {\rm O}({h^2}),\\
{}_x^{}D_L^\zeta u({x_i}) = {h^{ - \zeta }}\sum\limits_{k = 0}^{m - i + 1} {\vartheta _k^{(\zeta )}u({x_{i + k - 1}}} ) + {\rm O}({h^2}),
\end{array}
\end{equation}
where
\begin{equation}\label{eq56}
\vartheta _0^{(\zeta )} = \frac{\zeta }{2}\omega _0^{(\zeta )},\,\,\,\,\,\,\,\vartheta _k^{(\zeta )} = \frac{\zeta }{2}\omega _k^{(\zeta )} + \frac{{2 - \zeta }}{2}\omega _{k - 1}^{(\zeta )},\,\,\,\,\,k \ge 1.
\end{equation}
{\bf Lemma 5.} (\cite{tian2015class}) Suppose that $1 < \zeta  \le 2$, then the coefficients $\vartheta _k^{(\zeta )}$ satisfy
\begin{equation}\label{eq57}
\begin{array}{l}
\vartheta _0^{(\zeta )} = \frac{\zeta }{2} > 0,\,\,\,\,\vartheta _1^{(\zeta )} = \frac{{2 - \zeta  - {\zeta ^2}}}{2} < 0,\,\,\,\,\,\vartheta _2^{(\zeta )} = \frac{{\zeta ({\zeta ^2} + \zeta  - 4)}}{4},\\
1 \ge \vartheta _0^{(\zeta )} \ge \vartheta _3^{(\zeta )} \ge \vartheta _4^{(\zeta )} \ge ... \ge 0,\\
\sum\limits_{k = 0}^\infty  {\vartheta _k^{(\zeta )}}  = 0,\,\,\,\,\,\,\sum\limits_{k = 0}^m {\vartheta _k^{(\zeta )}}  < 0,\,\,\,\,\,\,m \ge 2.
\end{array}
\end{equation}

\section{Numerical method for the RFADE with delay}
In this section, we approximate the Riesz space fractional derivative and derive the Crank-Nicolson scheme of the equation. We partition the interval $[0, L]$ into an uniform mesh with the space step $h = L/M$ and the time step
$t = T/N$, where $M$, $N$ being two positive integers. The set of grid points are denoted by $x_{i} = ih$ and $t_{j}= j\kappa$ for $i=1,..., M$ and $j=1,..., M$.\\
We consider Eq (\ref{eq4}) and suppose $u({x_i},{t_j}) = u_i^j,\,\,\,\,f({x_i},{t_j}) = f_i^j$ and $\frac{{{\partial ^\gamma }g(x,t - \tau )}}{{\partial {t^\gamma }}} = \overline {g(x,t)}$. We take ${\eta _\alpha } = \frac{{{c_\alpha }\kappa {K_\alpha }{h^{ - \alpha }}}}{2}$ and ${\eta _\beta } = \frac{{{c_\beta }\kappa {K_\beta }{h^{ - \beta }}}}{2}$. Therefore, we have
\begin{equation}\label{eq59}
\begin{array}{*{20}{l}}
{u_i^j + {\eta _\alpha }(\sum\limits_{k = 0}^i {\varpi _k^{(\alpha )}} u_{i - k}^j + \sum\limits_{k = 0}^{M - i} {\varpi _k^{(\alpha )}} u_{i + k}^j) + {\eta _\beta }(\sum\limits_{k = 0}^{i + 1} {\vartheta _k^{(\beta )}} u_{i - k + 1}^j + \sum\limits_{k = 0}^{M - i + 1} {\vartheta _k^{(\beta )}} u_{i + k - 1}^j)}\\
{ = u_i^{j - 1} - {\eta _\alpha }(\sum\limits_{k = 0}^i {\varpi _k^{(\alpha )}} u_{i - k}^{j - 1} + \sum\limits_{k = 0}^{M - i} {\varpi _k^{(\alpha )}} u_{i + k}^{j - 1}) - {\eta _\beta }(\sum\limits_{k = 0}^{i + 1} {\vartheta _k^{(\beta )}} u_{i - k + 1}^{j - 1} + \sum\limits_{k = 0}^{M - i + 1} {\vartheta _k^{(\beta )}} u_{i + k - 1}^{j - 1})}\\
{ + \frac{\kappa }{2}(f_i^j + f_i^{j - 1}) - \frac{\kappa }{2}(\overline g _i^j + \overline g _i^{j - 1}),}
\end{array}
\end{equation}
we take
\begin{equation}\label{eq60}
\begin{array}{l}
A = \left( {\begin{array}{*{20}{c}}
{\varpi _0^{(\alpha )}}&0&0&{\begin{array}{*{20}{c}}
 \cdots &0
\end{array}}\\
{\varpi _1^{(\alpha )}}&{\varpi _0^{(\alpha )}}&0&{\begin{array}{*{20}{c}}
 \cdots &0
\end{array}}\\
 \vdots & \vdots & \vdots &{\begin{array}{*{20}{c}}
 \vdots & \vdots
\end{array}}\\
{\,\,\varpi _{M - 2}^{(\alpha )}}&{\varpi _{M - 3}^{(\alpha )}}& \cdots &{\,\,\,\, \cdots \varpi _0^{(\alpha )}}
\end{array}} \right),\\
B = \left( {\begin{array}{*{20}{c}}
{\vartheta _1^{(\beta )}}&{\vartheta _0^{(\beta )}}&0&{\begin{array}{*{20}{c}}
{\begin{array}{*{20}{c}}
0& \cdots
\end{array}}&0
\end{array}}\\
{\vartheta _2^{(\beta )}}&{\vartheta _1^{(\beta )}}&{\vartheta _0^{(\beta )}}&{\begin{array}{*{20}{c}}
0&{\begin{array}{*{20}{c}}
 \cdots &0
\end{array}}
\end{array}}\\
 \vdots & \vdots & \vdots &{\begin{array}{*{20}{c}}
 \vdots & \vdots
\end{array}}\\
{\,\,\vartheta _{M - 1}^{(\beta )}}&{\vartheta _{M - 2}^{(\beta )}}& \cdots &{\,\,\,\begin{array}{*{20}{c}}
 \cdots &{\vartheta _1^{(\beta )}}
\end{array}}
\end{array}} \right)
\end{array}
\end{equation}
and
\begin{equation}\label{eq61}
\begin{array}{l}
D = {\eta _\alpha }(A + {A^T}) + {\eta _\beta}(B + {B^T}),\\
{U^j} = {[u_1^j,{\mkern 1mu} u_2^j,{\mkern 1mu} ...,{\mkern 1mu} u_{M - 1}^j]^T}.
\end{array}
\end{equation}
Therefore, Eq. (\ref{eq59}) can be written as
\begin{equation}\label{eq62}
(I + D){U^j} = (I - D){U^{j - 1}} + {Q^j},
\end{equation}
where
\[{Q^j} = \left[ \begin{array}{l}
\frac{\kappa }{2}(f_1^j + f_1^{j - 1}) - \frac{\kappa }{2}(\overline g _1^j + \overline g _1^{j - 1}) - \Lambda \\
\frac{\kappa }{2}(f_2^j + f_2^{j - 1}) - \frac{\kappa }{2}(\overline g _2^j + \overline g _2^{j - 1}) - {\Lambda _2}\\
\frac{\kappa }{2}(f_3^j + f_3^{j - 1}) - \frac{\kappa }{2}(\overline g _3^j + \overline g _3^{j - 1}) - {\Lambda _3}\\
 \vdots \\
\frac{\kappa }{2}(f_{M - 1}^j + f_{M - 1}^{j - 1}) - \frac{\kappa }{2}(\overline g _{M - 1}^j + \overline g _{M - 1}^{j - 1}) - {\Lambda _{M - 1}}
\end{array} \right]\]
and
\[{\Lambda _s} = ({\eta _\alpha }\varpi _s^{(\alpha )} + {\eta _\beta }\vartheta _{s + 1}^{(\beta )})(u_0^j + u_0^{j - 1}) + ({\eta _\alpha }\varpi _{M - s}^{(\alpha )} + {\eta _\beta }\vartheta _{M - s + 1}^{(\beta )})(u_M^j + u_M^{j - 1}),\,\,\,\,s = 2,3,...,M - 1,\]

\[\Lambda  = ({\eta _\alpha }\varpi _1^{(\alpha )} + {\eta _\beta }(\vartheta _0^{(\beta )} + \vartheta _2^{(\beta )}))(u_0^j + u_0^{j - 1}) + ({\eta _\alpha }\varpi _{M - 1}^{(\alpha )} + {\eta _\beta }\vartheta _M^{(\beta )})(u_M^j + u_M^{j - 1}).\]

\subsection{Stability of method}
{\bf Lemma 6.} (\cite{thomas2013numerical}) Let $A$ be an $m - 1$ order positive define matrix, then for any parameter $\nu  \ge 0$, the following two inequalities hold
\[{\left\| {{{(I + \nu A)}^{ - 1}}} \right\|_\infty } \le 1,\,\,\,\,\,\,{\left\| {{{(I + \nu A)}^{ - 1}}(I - \nu A)} \right\|_\infty } \le 1\]

{\bf Theorem 7.} Matrix $D$ (\ref{eq61}) for $h < {( - \frac{{3\beta (\beta  - 1)(2 - \beta )(3 + \beta )\cos (\frac{{\alpha \pi }}{2})}}{{4{{(\frac{3}{2})}^\alpha }\alpha (8\alpha  - 5)\cos (\frac{{\beta \pi }}{2})}})^{\frac{1}{{\beta  - \alpha }}}}$, is strictly diagonally dominant matrix.\\
{\bf Proof.} We have
\[{D_{i,j}} = \left\{ \begin{array}{l}
{\eta _\alpha }\varpi _{j - i}^{(\alpha )} + {\eta _\beta }\vartheta _{j - i + 1}^{(\beta )},\,\,\,\,\,\,\,\,\,\,\,\,\,\,\,\,\,\,\,\,\,\,\,\,j > i + 1,\\
{\eta _\alpha }\varpi _1^{(\alpha )} + {\eta _\beta }(\vartheta _0^{(\beta )} + \vartheta _2^{(\beta )}),\,\,\,\,\,\,\,\,\,\,j = i + 1,\\
2{\eta _\alpha }\varpi _0^{(\alpha )} + 2{\eta _\beta }\vartheta _1^{(\beta )},\,\,\,\,\,\,\,\,\,\,\,\,\,\,\,\,\,\,\,\,\,\,\,\,\,\,\,j = i,\\
{\eta _\alpha }\varpi _1^{(\alpha )} + {\eta _\beta }(\vartheta _0^{(\beta )} + \vartheta _2^{(\beta )}),\,\,\,\,\,\,\,\,\,\,j = i - 1,\\
{\eta _\alpha }\varpi _{i - j}^{(\alpha )} + {\eta _\beta }\vartheta _{i - j + 1}^{(\beta )}.\,\,\,\,\,\,\,\,\,\,\,\,\,\,\,\,\,\,\,\,\,\,\,\,j < i + 1,
\end{array} \right.\]
where ${\eta _\alpha }>0$ for $0<\alpha<1$ and ${\eta _\beta }<0$ for $1 < \beta  \le 2$. According to Lemma 3, when $k \ge 4$, ${\eta _\alpha }\varpi _k^{(\alpha )} < 0$ and according to Lemma 5, when $k \ge 4$,  ${\eta _\beta }\vartheta _k^{(\beta )} < 0$. For $k=2,3$, if $0 < \alpha   \le \frac{5}{8},\,\,1 < \beta  < 2$, then ${\eta _\alpha }\varpi _2^{(\alpha )} < 0$ and ${\eta _\alpha }\varpi _3^{(\alpha )} < 0$. Therefore, ${D_{i,j}} < 0$ when
$j > i + 1$ or $j < i - 1$. According to the Lemmas 3 and 5, we have $\varpi _0^{(\alpha )} > 0,\,\,\,\vartheta _1^{(\beta )} < 0$,
 thus
\[2{\eta _\alpha }\varpi _0^{(\alpha )} + 2{\eta _\beta }\vartheta _1^{(\beta )} > 0.\]
Therefore ${D_{i,i}} > 0$. For ${D_{i,i+1}}$ and ${D_{i,i-1}}$, we have
\[\begin{array}{l}
\vartheta _0^{(\beta )} + \vartheta _2^{(\beta )} = \frac{\beta }{2} + \frac{{\beta ({\beta ^2} + \beta  - 4)}}{4} = \frac{{\beta (\beta  + 2)(\beta  - 1)}}{4} > 0,\\
\varpi _1^{(\alpha )} < 0.
\end{array}\]
Since ${\eta _\alpha }>0$  and ${\eta _\beta }<0$, then
\[{D_{i,i + 1}} = {D_{i,i - 1}} = {\eta _\alpha }\varpi _1^{(\alpha )} + {\eta _\beta }(\vartheta _0^{(\beta )} + \vartheta _2^{(\beta )}) < 0.\]
For a given $i$, we can write
\begin{equation}\label{eq4000}
\begin{array}{l}
\sum\limits_{j = 1,\,j \ne i}^{M - 1} {\left| {{D_{i,j}}} \right| = } \sum\limits_{j = 1}^{i - 2} {\left| {{D_{i,j}}} \right|}  + \sum\limits_{j = i + 2}^{M - 1} {\left| {{D_{i,j}}} \right|}  + \left| {{D_{i,i - 1}}} \right| + \left| {{D_{i,i + 1}}} \right|\\
\,\,\,\,\,\,\,\,\,\,\,\,\,\,\,\,\,\,\,\,\,\,\,\,\, =  - \sum\limits_{j = 1}^{i - 2} {({\eta _\alpha }\varpi _{i - j}^{(\alpha )} + {\eta _\beta }\vartheta _{i - j + 1}^{(\beta )}) - } \sum\limits_{j = i + 2}^{M - 1} {({\eta _\alpha }\varpi _{j - i}^{(\alpha )} + {\eta _\beta }\vartheta _{j - i + 1}^{(\beta )})} \\
\,\,\,\,\,\,\,\,\,\,\,\,\,\,\,\,\,\,\,\,\,\,\,\,\,\,\,\,\, - 2{\eta _\alpha }\varpi _1^{(\alpha )} - 2{\eta _\beta }(\vartheta _0^{(\beta )} + \vartheta _2^{(\beta )})\\
\,\,\,\,\,\,\,\,\,\,\,\,\,\,\,\,\,\,\,\,\,\,\,\,\,\, < \sum\limits_{j =  - \infty }^{i - 2} {({\eta _\alpha }\varpi _{i - j}^{(\alpha )} + {\eta _\beta }\vartheta _{i - j + 1}^{(\beta )}) - \sum\limits_{j = i + 2}^{ + \infty } {({\eta _\alpha }\varpi _{j - i}^{(\alpha )} + {\eta _\beta }\vartheta _{j - i + 1}^{(\beta )})} } \\
\,\,\,\,\,\,\,\,\,\,\,\,\,\,\,\,\,\,\,\,\,\,\,\,\,\,\,\,\,\, - 2{\eta _\alpha }\varpi _1^{(\alpha )} - 2{\eta _\beta }(\vartheta _0^{(\beta )} + \vartheta _2^{(\beta )})\\
\,\,\,\,\,\,\,\,\,\,\,\,\,\,\,\,\,\,\,\,\,\,\,\,\,\,\, =  - 2{\eta _\alpha }\sum\limits_{k = 2}^{ + \infty } {\varpi _k^{(\alpha )}}  - 2{\eta _\beta }\sum\limits_{k = 3}^{ + \infty } {\vartheta _k^{(\beta )} - } 2{\eta _\alpha }\varpi _1^{(\alpha )} - 2{\eta _\beta }(\vartheta _0^{(\beta )} + \vartheta _2^{(\beta )})\\
\,\,\,\,\,\,\,\,\,\,\,\,\,\,\,\,\,\,\,\,\,\,\,\,\,\,\, =  - 2{\eta _\alpha }\sum\limits_{k = 0}^{ + \infty } {\varpi _k^{(\alpha )} - } 2{\eta _\beta }\sum\limits_{k = 0}^{ + \infty } {\vartheta _k^{(\beta )} + 2{\eta _\alpha }\varpi _1^{(\alpha )} + } 2{\eta _\beta }\vartheta _1^{(\beta )}\\
\,\,\,\,\,\,\,\,\,\,\,\,\,\,\,\,\,\,\,\,\,\,\,\,\,\,\, = 2{\eta _\alpha }\varpi _1^{(\alpha )} + 2{\eta _\beta }\vartheta _1^{(\beta )} = \left| {{D_{i,i}}} \right|,
\end{array}
\end{equation}
Therefore
\[\sum\limits_{j = 1,\,j \ne i}^{M - 1} {\left| {{D_{i,j}}} \right| < } \left| {{D_{i,i}}} \right|.\]
Also, if $\frac{5}{8} < \alpha  < 1,\,\,1 < \beta  < 2$, for $h < {( - \frac{{3\beta (\beta  - 1)(2 - \beta )(3 + \beta )\cos (\frac{{\alpha \pi }}{2})}}{{4{{(\frac{3}{2})}^\alpha }\alpha (8\alpha  - 5)\cos (\frac{{\beta \pi }}{2})}})^{\frac{1}{{\beta  - \alpha }}}}$, we can write
\[\begin{array}{l}
{\eta _\alpha }\varpi _2^{(\alpha )} + {\eta _\beta }\vartheta _3^{(\beta )} < 0,\\
{\eta _\alpha }\varpi _3^{(\alpha )} + {\eta _\beta }\vartheta _4^{(\beta )} < 0.
\end{array}\]
Therefore, ${D_{i,j}} < 0$ when $j > i + 1$ or $j < i - 1$. Then relation (\ref{eq4000}) is valid for $\frac{5}{8}< \alpha  < 1,\,\,1 < \beta  < 2$. Thus, matrix $D$ is strictly diagonally dominant matrix.\\
{\bf Lemma 7.} The matrix $D$ (\ref{eq61}) for $h < {( - \frac{{3\beta (\beta  - 1)(2 - \beta )(3 + \beta )\cos (\frac{{\alpha \pi }}{2})}}{{4{{(\frac{3}{2})}^\alpha }\alpha (8\alpha  - 5)\cos (\frac{{\beta \pi }}{2})}})^{\frac{1}{{\beta  - \alpha }}}}$, is symmetric positive definite.\\
{\bf Proof.} The matrix $D$ by using Eq. \ref{eq61} is clearly symmetric. Let ${\eta _0}$ be one eigenvalue of the matrix $D$.  By using the (\cite{varga2010gervsgorin}), we can write
\[\left| {{\eta _0} - {D_{i,i}}} \right| \le \sum\limits_{j = 1,\,j \ne i}^{M - 1} {\left| {{D_{i,j}}} \right|,} \]
or
\[{D_{i,i}} - \sum\limits_{j = 1,\,j \ne i}^{M - 1} {\left| {{D_{i,j}}} \right|}  \le {\eta _0} \le {D_{i,i}} +\sum\limits_{j = 1,\,j \ne i}^{M - 1} {\left| {{D_{i,j}}} \right|} ,\]
by using Theorem 7, we have
\[{\eta _0} \ge {D_{i,i}} - \sum\limits_{j = 1,\,j \ne i}^{M - 1} {\left| {{D_{i,j}}} \right|}  \ge 0,\]
thus D is positive definite.\\
{\bf Theorem 8.} The numerical method (\ref{eq62}) for $h < {( - \frac{{3\beta (\beta  - 1)(2 - \beta )(3 + \beta )\cos (\frac{{\alpha \pi }}{2})}}{{4{{(\frac{3}{2})}^\alpha }\alpha (8\alpha  - 5)\cos (\frac{{\beta \pi }}{2})}})^{\frac{1}{{\beta  - \alpha }}}}$, is stable.\\
{\bf Proof.} Let ${U^j}$ be the numerical solution and ${u^j}$ be the exact solution. Since the matrix (I + D) for $h < {( - \frac{{3\beta (\beta  - 1)(2 - \beta )(3 + \beta )\cos (\frac{{\alpha \pi }}{2})}}{{4{{(\frac{3}{2})}^\alpha }\alpha (8\alpha  - 5)\cos (\frac{{\beta \pi }}{2})}})^{\frac{1}{{\beta  - \alpha }}}}$, is invertible, then if we take
\begin{equation}\label{eq63}
{\varepsilon ^j} = {U^j} - {u^j},\,\,P = {(I + D)^{ - 1}}{(I-D)}.
\end{equation}
Therefore, we have
\[{\varepsilon ^j} = P{\varepsilon ^{j - 1}}.\]
Then
\[{\varepsilon ^j} = {P^{j - 1}}{\varepsilon ^0}.\]
By using Lemma 6, we have
\[{\left\| {{\varepsilon ^j}} \right\|_\infty } = {\left\| {{P^{j - 1}}{\varepsilon ^0}} \right\|_\infty } \le \left\| P \right\|_\infty ^{j - 1}{\left\| {{\varepsilon ^0}} \right\|_\infty } = \left\| {{{(I + D)}^{ - 1}}(I-D)} \right\|_\infty ^{j - 1}{\left\| {{\varepsilon ^0}} \right\|_\infty } \le {\left\| {{\varepsilon ^0}} \right\|_\infty }.\]
Therefore, the numerical method (\ref{eq62}) is stable.\\

\subsection{Convergence of method}
In this subsection, we study the local truncation error of method, therefore we can write
\begin{equation}\label{eq73}
\begin{array}{l}
\dfrac{{u({x_i},{t_j}) - u({x_i},{t_{j - 1}})}}{\kappa } = \dfrac{1}{2}({K_\alpha }\dfrac{{{\partial ^\alpha }u({x_i},{t_j})}}{{\partial {{\left| x \right|}^\alpha }}} + {K_\beta }\dfrac{{{\partial ^\beta }u({x_i},{t_j})}}{{\partial {{\left| x \right|}^\beta }}})\\
\,\,\,\,\,\,\,\,\,\,\,\,\,\,\,\,\,\,\,\,\,\,\,\,\,\,\,\,\,\,\,\,\,\,\,\,\,\,\,\,\,\,\,\,\,\, + \dfrac{1}{2}({K_\alpha }\dfrac{{{\partial ^\alpha }u({x_i},{t_{j - 1}})}}{{\partial {{\left| x \right|}^\alpha }}} + {K_\beta }\dfrac{{{\partial ^\beta }u({x_i},{t_{j - 1}})}}{{\partial {{\left| x \right|}^\beta }}}) + {\rm O}({\kappa ^2})\\
{K_\alpha }\dfrac{{{\partial ^\alpha }u({x_i},{t_j})}}{{\partial {{\left| x \right|}^\alpha }}} = {\eta _\alpha }(\sum\limits_{k = 0}^i {\varpi _k^{(\alpha )}} u_{i - k}^j + \sum\limits_{k = 0}^{M - i} {\varpi _k^{(\alpha )}} u_{i + k}^j) + {\rm O}({h^2}),\\
{K_\beta }\dfrac{{{\partial ^\beta }u({x_i},{t_j})}}{{\partial {{\left| x \right|}^\beta }}} = {\eta _\beta }(\sum\limits_{k = 0}^{i + 1} {\vartheta _k^{(\beta )}} u_{i - k + 1}^j + \sum\limits_{k = 0}^{M - i + 1} {\vartheta _k^{(\beta )}} u_{i + k - 1}^j) + {\rm O}({h^2}).
\end{array}
\end{equation}
Thus the local truncation error of (\ref{eq62}), will be as the following form
\[{T_{i,j}} = {\rm O}({\kappa ^3} + \kappa {h^2}).\]
{\bf Theorem 9.} Let ${U^j}$ be the numerical solution and ${u^j}$ be the exact solution of (\ref{eq59}), if $h < {( - \frac{{3\beta (\beta  - 1)(2 - \beta )(3 + \beta )\cos (\frac{{\alpha \pi }}{2})}}{{4{{(\frac{3}{2})}^\alpha }\alpha (8\alpha  - 5)\cos (\frac{{\beta \pi }}{2})}})^{\frac{1}{{\beta  - \alpha }}}}$, Then we have
\begin{equation}\label{eq74}
{\left\| {{U^j} - {u^j}} \right\|_\infty } \le C{\rm O}({\kappa ^2}  + {h^2}),
\end{equation}
where $C$ denotes a positive constant.\\
{\bf Proof.} We consider the Eq (\ref{eq59}), thus we can write
\begin{equation}\label{eq75}
\begin{array}{*{20}{l}}
{u_i^j + {\eta _\alpha }(\sum\limits_{k = 0}^i {\varpi _k^{(\alpha )}} u_{i - k}^j + \sum\limits_{k = 0}^{M - i} {\varpi _k^{(\alpha )}} u_{i + k}^j) + {\eta _\beta }(\sum\limits_{k = 0}^{i + 1} {\vartheta _k^{(\beta )}} u_{i - k + 1}^j + \sum\limits_{k = 0}^{M - i + 1} {\vartheta _k^{(\beta )}} u_{i + k - 1}^j)}\\
{ = u_i^{j - 1} - {\eta _\alpha }(\sum\limits_{k = 0}^i {\varpi _k^{(\alpha )}} u_{i - k}^{j - 1} + \sum\limits_{k = 0}^{M - i} {\varpi _k^{(\alpha )}} u_{i + k}^{j - 1}) - {\eta _\beta }(\sum\limits_{k = 0}^{i + 1} {\vartheta _k^{(\beta )}} u_{i - k + 1}^{j - 1} + \sum\limits_{k = 0}^{M - i + 1} {\vartheta _k^{(\beta )}} u_{i + k - 1}^{j - 1})}\\
{ + \frac{\kappa }{2}(f_i^j + f_i^{j - 1}) - \frac{\kappa }{2}(\overline g _i^j + \overline g _i^{j - 1})}
\end{array}
\end{equation}
and
\begin{equation}\label{eq76}
\begin{array}{*{20}{l}}
{U_i^j + {\eta _\alpha }(\sum\limits_{k = 0}^i {\varpi _k^{(\alpha )}} U_{i - k}^j + \sum\limits_{k = 0}^{M - i} {\varpi _k^{(\alpha )}} U_{i + k}^j) + {\eta _\beta }(\sum\limits_{k = 0}^{i + 1} {\vartheta _k^{(\beta )}} U_{i - k + 1}^j + \sum\limits_{k = 0}^{M - i + 1} {\vartheta _k^{(\beta )}} U_{i + k - 1}^j)}\\
{ = U_i^{j - 1} - {\eta _\alpha }(\sum\limits_{k = 0}^i {\varpi _k^{(\alpha )}} U_{i - k}^{j - 1} + \sum\limits_{k = 0}^{M - i} {\varpi _k^{(\alpha )}} U_{i + k}^{j - 1}) - {\eta _\beta }(\sum\limits_{k = 0}^{i + 1} {\vartheta _k^{(\beta )}} U_{i - k + 1}^{j - 1} + \sum\limits_{k = 0}^{M - i + 1} {\vartheta _k^{(\beta )}} U_{i + k - 1}^{j - 1})}\\
{ + \frac{\kappa }{2}(f_i^j + f_i^{j - 1}) - \frac{\kappa }{2}(\overline g _i^j + \overline g _i^{j - 1}).}
\end{array}
\end{equation}
Let $e_i^j = U_i^j - u_i^j$ and by using (\ref{eq75}) and (\ref{eq76}), we have
\begin{equation}\label{eq77}
\begin{array}{l}
e_i^j + {\eta _\alpha }(\sum\limits_{k = 0}^i {\varpi _k^{(\alpha )}} e_{i - k}^j + \sum\limits_{k = 0}^{M - i} {\varpi _k^{(\alpha )}} e_{i + k}^j) + {\eta _\beta }(\sum\limits_{k = 0}^{i + 1} {\vartheta _k^{(\beta )}} e_{i - k + 1}^j + \sum\limits_{k = 0}^{M - i + 1} {\vartheta _k^{(\beta )}} e_{i + k - 1}^j)\\
 = e_i^{j - 1} - {\eta _\alpha }(\sum\limits_{k = 0}^i {\varpi _k^{(\alpha )}} e_{i - k}^{j - 1} + \sum\limits_{k = 0}^{M - i} {\varpi _k^{(\alpha )}} e_{i + k}^{j - 1}) - {\eta _\beta }(\sum\limits_{k = 0}^{i + 1} {\vartheta _k^{(\beta )}} e_{i - k + 1}^{j - 1} + \sum\limits_{k = 0}^{M - i + 1} {\vartheta _k^{(\beta )}} e_{i + k - 1}^{j - 1}).
\end{array}
\end{equation}
We can write
\[(I + D){\varepsilon ^j} = (I -D){\varepsilon ^{j-1}} + {\rm O}({\kappa ^3} + \kappa {h^2})\chi ,\]
where
\[{\varepsilon ^j} = {[e_1^j,e_2^j,...,e_n^j]^T},{\mkern 1mu} {\mkern 1mu} {\mkern 1mu} \chi  = {[1,1,...,1]^T},{\mkern 1mu} {\mkern 1mu} {\rm{ }}D = {\eta _\alpha }(A + {A^T}) + {\eta _\beta }(B + {B^T}).\]
Let we take
\[L = {(I + D)^{ - 1}(I+D)},\,\,X= {\rm O}({\kappa ^3} + \kappa {h^2}){(I + D)^{ - 1}}.\]
So, we can rewrite
\[{\varepsilon ^j} = L{\varepsilon ^{j - 1}} + X,\]
with iterating and by using initial condition of main problem, we have
\[{\varepsilon ^j} = ({L^{j - 1}} + {L^{j - 2}} + ... + I)X.\]
By using Lemma 6, we have
\[\begin{array}{l}
{\left\| {{\varepsilon ^j}} \right\|_\infty } \le ({\left\| {{L^{j - 1}}} \right\|_\infty } + {\left\| {{L^{j - 2}}} \right\|_\infty } + ... + {\left\| I \right\|_\infty }){\left\| X \right\|_\infty }\\
\,\,\,\,\,\,\,\,\,\,\,\,\, = (\left\| L \right\|_\infty ^{j - 1} + \left\| L \right\|_\infty ^{j - 2} + ... + {\left\| I \right\|_\infty }){\left\| X \right\|_\infty }\\
\,\,\,\,\,\,\,\,\,\,\,\,\, \le (1 + 1 + ... + 1){\left\| X \right\|_\infty }\\
\,\,\,\,\,\,\,\,\,\,\,\,\, \le j{\rm O}({\kappa ^3} + \kappa {h^2}) = T{\rm O}({\kappa ^2}  + {h^2}).
\end{array}\]
Therefore, we have
\[{\left\| {{\varepsilon ^j}} \right\|_\infty } \le C{\rm O}({\kappa ^2}  + {h^2}).\]
\section{Test examples}
Since this model (\ref{eq4}) has not been considered so far, thus we consider three cases in this section. We consider the Riesz space fractional advection-dispersion equation without a source term and delay (${K_\gamma } = 0,\,\,\tau  = 0,\,\,\,f(x,t) = 0$). Next, the Riesz space fractional advection-dispersion equation with a source term and without delay (${K_\gamma } = 0,\,\,\tau  = 0$) is considered. Finally, an example with delay and a source term is given to demonstrate the applicability and accuracy of the presented numerical method.
\\
{\bf Example 1.} (\cite{yang2010numerical}) We consider the following Riesz space fractional advection-dispersion
equation
\begin{equation}\label{eq000001}
\begin{array}{l}
\dfrac{{\partial u(x,t)}}{{\partial t}} = {K_\alpha }\dfrac{{{\partial ^\alpha }u(x,t)}}{{\partial {{\left| x \right|}^\alpha }}} + {K_\beta }\dfrac{{{\partial ^\beta }u(x,t)}}{{\partial {{\left| x \right|}^\beta }}},\,\,\,\,0 < x < \pi ,\,\,\,\,0 < t \le T,\\
u(x,0) = {x^2}(\pi  - x),{\mkern 1mu} {\mkern 1mu} \\
u(0,t) = u(\pi ,t) = 0.{\mkern 1mu} {\mkern 1mu} \,
\end{array}
\end{equation}
where, $0<\alpha<1,\,\,1<\beta \le 2$. From \cite{yang2010numerical} we know the exact solution of  \ref{eq000001}  is given by
\begin{equation}\label{eq000002}
u(x,t) = \sum\limits_{n = 1}^\infty  {\left[ {\frac{8}{{{n^3}}}{{( - 1)}^{n + 1}} - \frac{4}{{{n^3}}}} \right]\sin (nx)\exp ( - \left[ {{K_\alpha }{{({n^2})}^{\frac{\alpha }{2}}} + {K_\beta }{{({n^2})}^{\frac{\beta }{2}}}} \right]t)} .
\end{equation}
The absolute errors of MTM scheme given in Yang et al.\cite{yang2010numerical} and presented scheme are shown in Tables 1 and they are compared for different values of $h$ and $\alpha  = 0.4,\,\,\beta  = 1.8,\,\,\,t = 0.4$ and ${K_\alpha } = {K_\beta } = 0.25$.
\begin{figure}[h!]
	\begin{flushleft}
		\includegraphics[width=13cm]{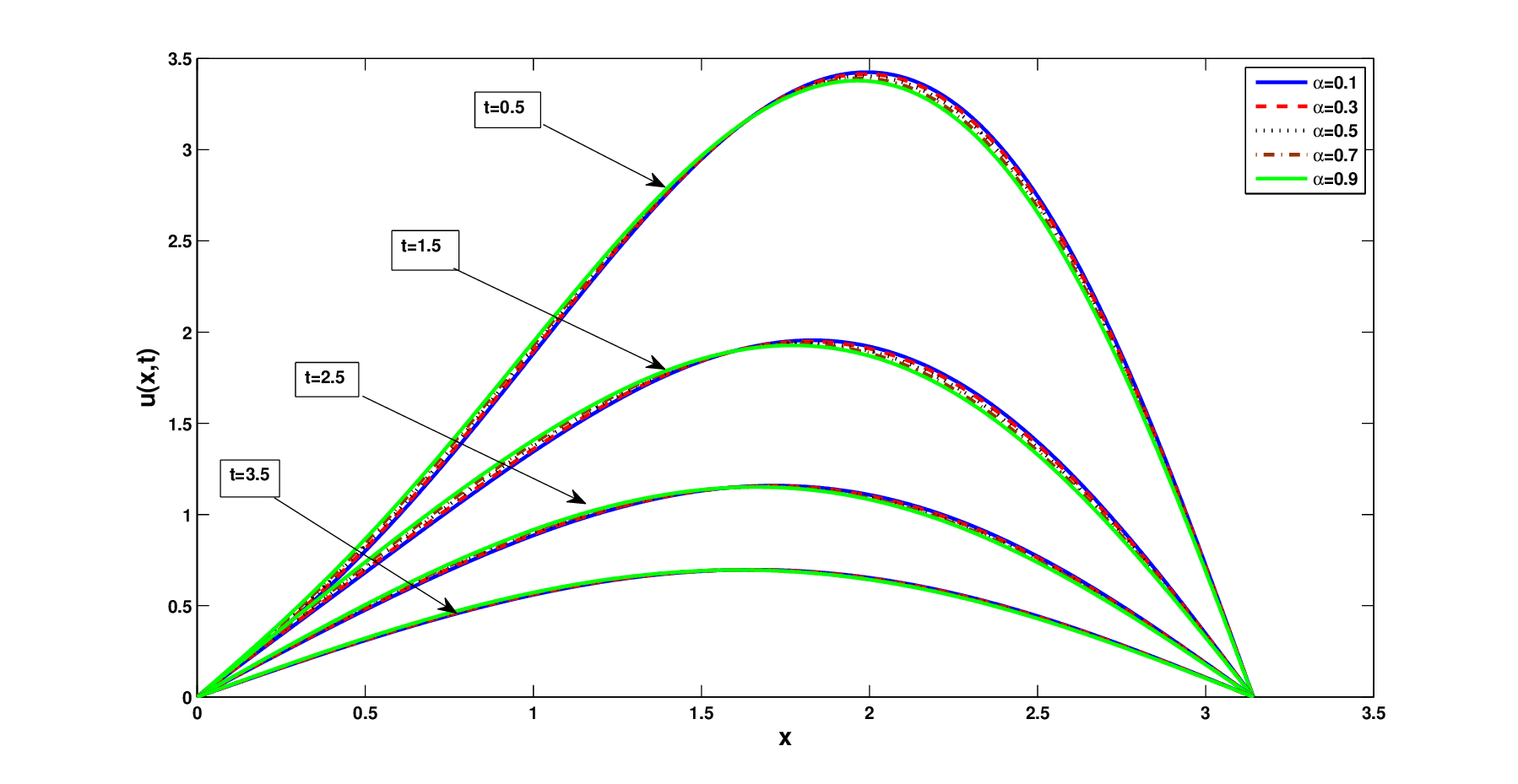}
		\caption{{\small The numerical approximation by the presented method for example 1, with different values of $\alpha = 0.1, 0.3, 0.5, 0.7, 0.9$, when $\beta=1.8$ and $T=0.5, 1.5, 2.5, 3.5$.}
		}\label{5}
	\end{flushleft}
\end{figure}

\begin{figure}[h!]
	\begin{flushleft}
		\includegraphics[width=13cm]{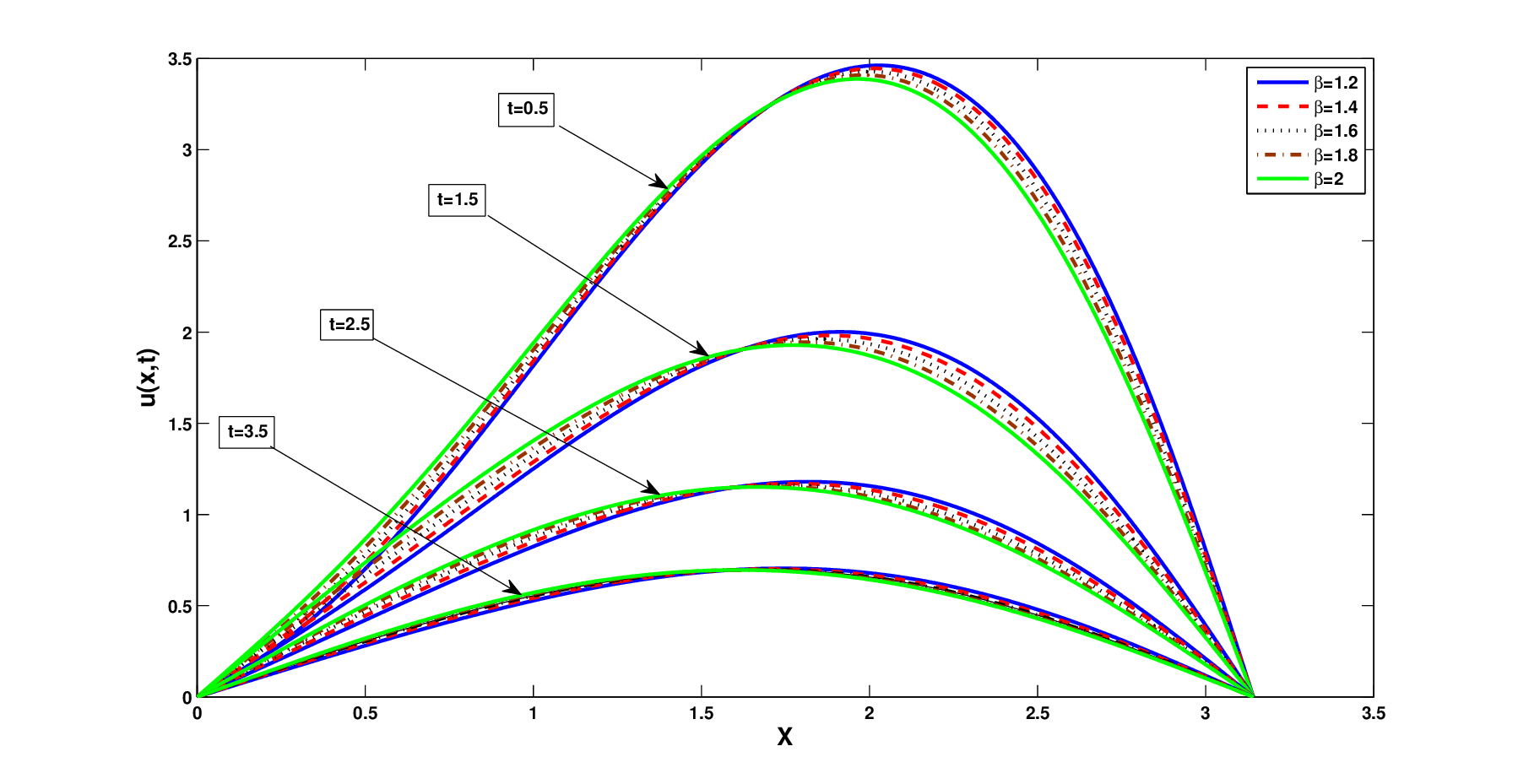}
		\caption{{\small The numerical approximation by the presented method for example 1, with different values of $\beta = 1.2, 1.4, 1.6, 1.8, 2$, when $\alpha=0.4$ and $T=0.5, 1.5, 2.5, 3.5$.}
		}\label{5}
	\end{flushleft}
\end{figure}

\begin{figure}[h!]
	\begin{flushleft}
		\includegraphics[width=13cm]{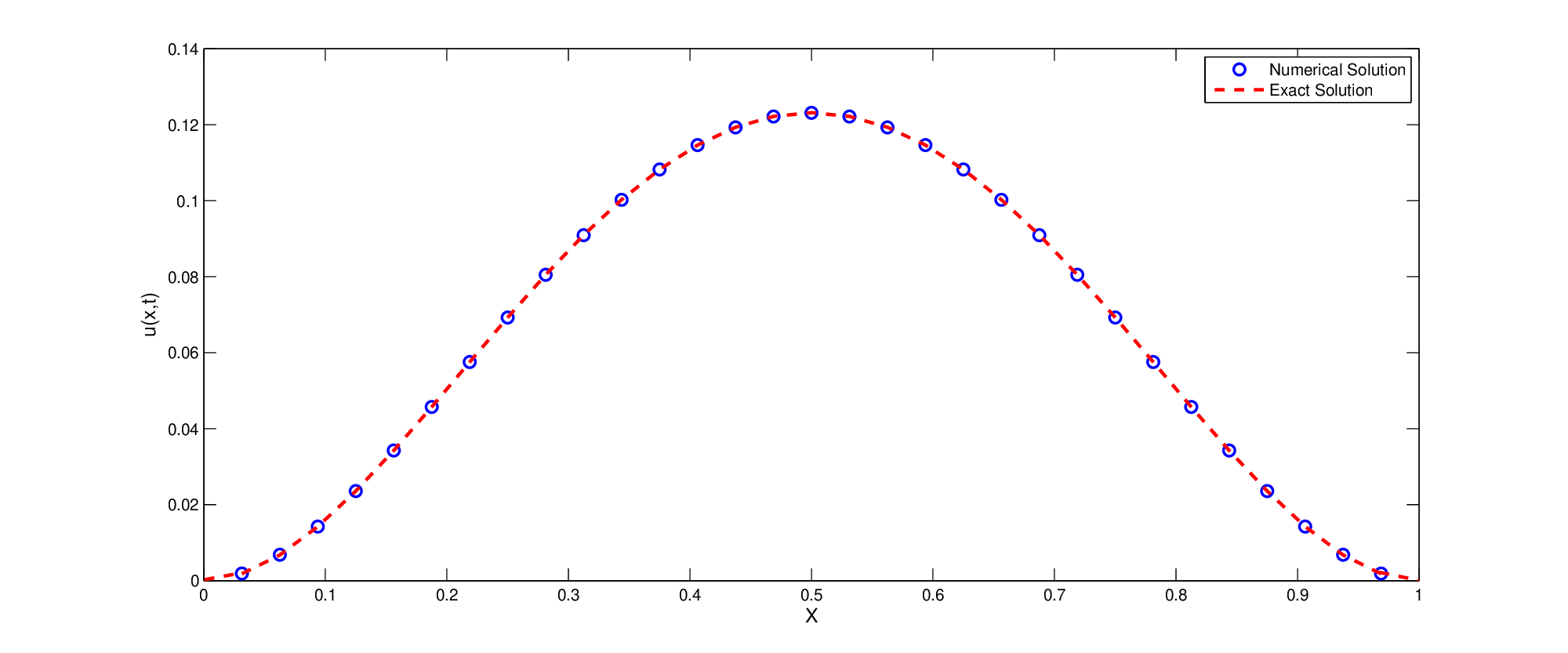}
		\caption{{\small The numerical approximation and exact solution by the presented method for example3, for $\alpha = 0.8$, $\gamma = 0.7$, $\beta = 1.2$, when T = 1.}
		}
	\end{flushleft}
\end{figure}

\begin{figure}[h!]
	\begin{flushleft}
		\includegraphics[width=13cm]{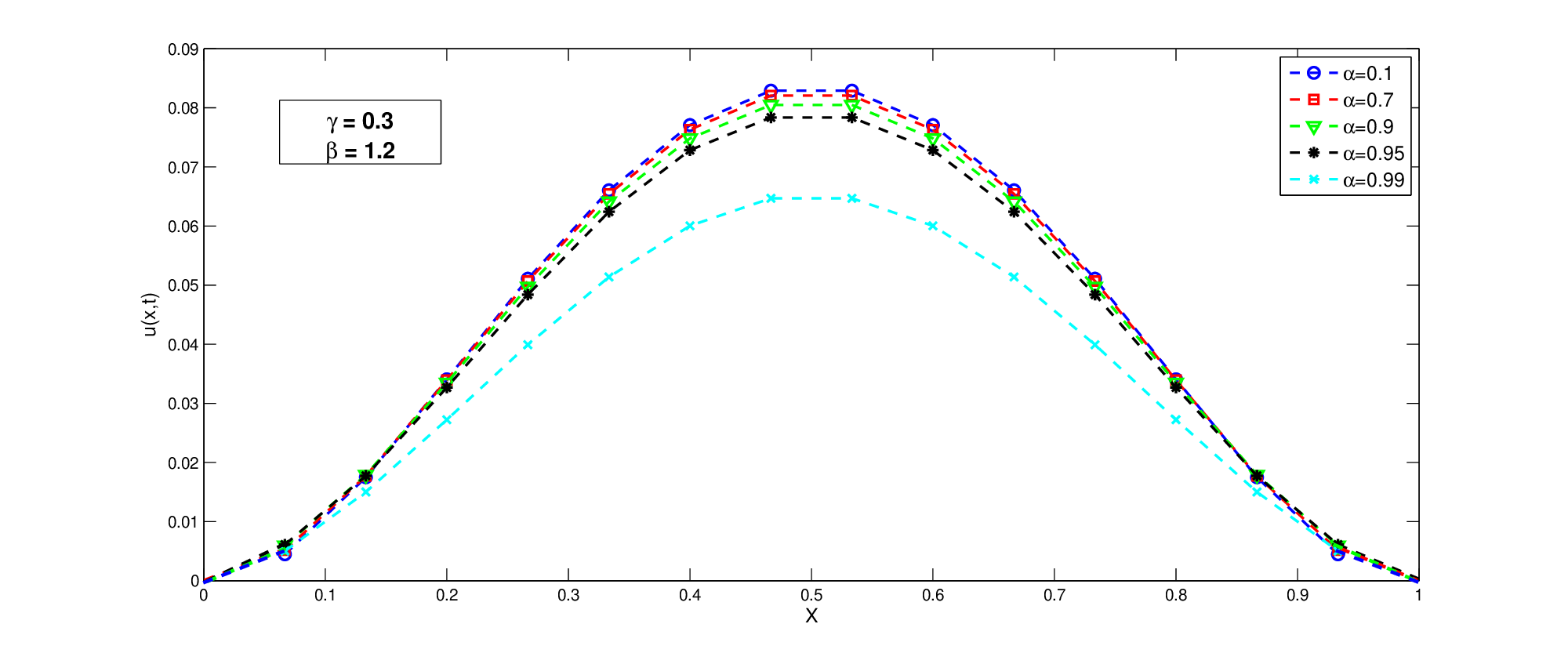}
		\caption{{\small The numerical approximation by the presented method for example3, for various $\alpha = 0.1, 0.7, 0.9, 0.95, 0.99 $, when T = 1 and $\gamma=0.3$, $\beta=1.2$.}
		}\label{5}
	\end{flushleft}
\end{figure}


\begin{figure}[h!]
	\begin{flushleft}
		\includegraphics[width=13cm]{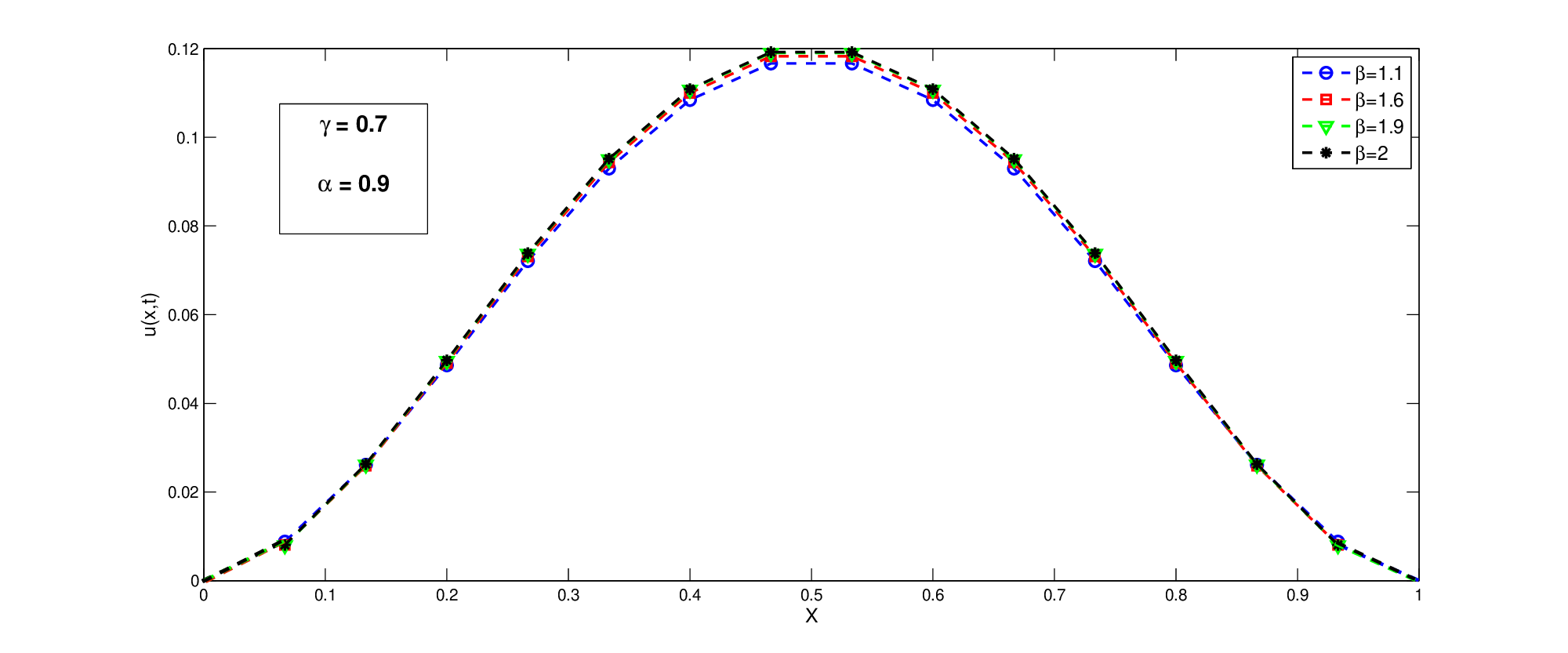}
		\caption{{\small The numerical approximation by the presented method for example3, for various $\beta = 1.1, 1.6, 1.9, 2 $, when T = 1 and $\gamma=0.7$, $\alpha=0.9$.}
		}\label{5}
	\end{flushleft}
\end{figure}


\begin{figure}[h!]
	\begin{flushleft}
		\includegraphics[width=13cm]{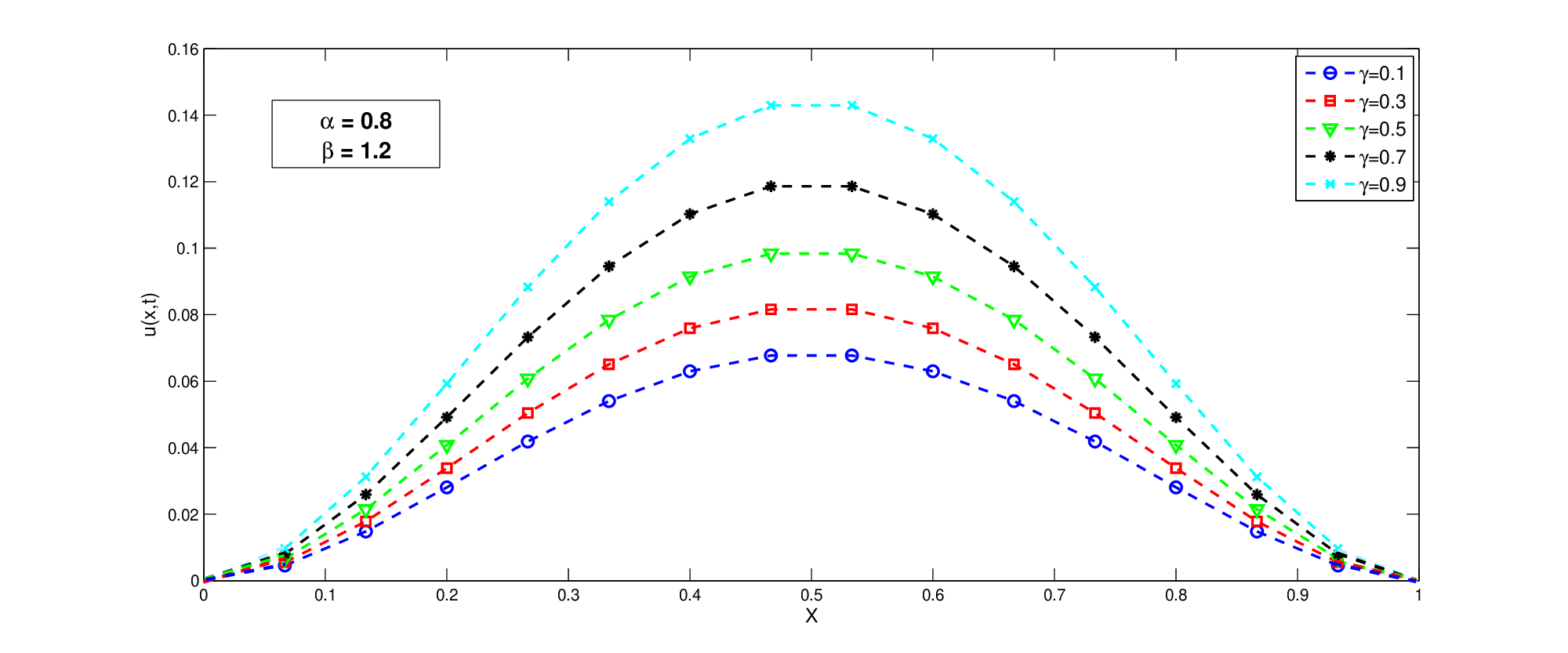}
		\caption{{\small The numerical approximation by the presented method for example3, for various $\gamma = 0.1, 0.3, 0.5, 0.7, 0.9 $, when T = 1 and $\alpha=0.8$, $\beta=1.2$.}
		}\label{5}
	\end{flushleft}
\end{figure}
From Table 1 and Figs 2-3, it can be observed that the errors of the presented method are significantly compared with the literature.
\begin{table}[ht]
\centering
\caption{The absolute errors of the present scheme and
numerical method in Yang et al.\cite{yang2010numerical} with $\alpha  = 0.4,\,\,\beta  = 1.8,\,\,\,t = 0.4$ for \ref{eq000001}.}
\begin{tabular}[t]{ccccc}
\hline
$h$&Proposed Method&order&MTM \cite{yang2010numerical}&order\\
\hline
$\frac{\pi }{{10}}$&1.55006e-003&-&1.9963e-002&-\\
$\frac{\pi }{{20}}$&3.9148e-004&1.9796&5.1683e-003&1.9496\\
$\frac{\pi }{{40}}$&9.8526e-005&1.9898&1.3386e-003&1.9490\\
$\frac{\pi }{{80}}$&2.4633e-005&1.9999&3.3519e-004&1.9977\\
\hline
\end{tabular}
\end{table}%

{\bf Example 2.} \cite{yang2010numerical} We consider the following Riesz space fractional advection-dispersion equation with a source term:
\begin{equation}\label{eq000002}
\begin{array}{l}
\dfrac{{\partial u(x,t)}}{{\partial t}} = {K_\alpha }\dfrac{{{\partial ^\alpha }u(x,t)}}{{\partial {{\left| x \right|}^\alpha }}} + {K_\beta }\dfrac{{{\partial ^\beta }u(x,t)}}{{\partial {{\left| x \right|}^\beta }}} + f(x,t),\,\,\,\,0 < t \le T,\,\,\,0 < x < 1,\\
u(x,0) = 0,{\mkern 1mu} {\mkern 1mu} \\
u(0,t) = u(1,t) = 0,{\mkern 1mu}
\end{array}
\end{equation}
where
\[\begin{array}{l}
f(x,t) = {t^{\alpha  - 1}}{e^{\beta t}}(\alpha  + \beta t){x^2}{(1 - x)^2}\\
\,\,\,\,\,\,\,\,\,\,\,\,\,\,\, + \dfrac{{{K_\alpha }{t^\alpha }{e^{\beta t}}}}{{2\cos (\dfrac{{\alpha \pi }}{2})}}\left[ {\dfrac{{24({{(1 - x)}^{4 - \alpha }} + {x^{4 - \alpha }})}}{{\Gamma (5 - \alpha )}} - \dfrac{{12({{(1 - x)}^{3 - \alpha }} + {x^{3 - \alpha }})}}{{\Gamma (4 - \alpha )}} + \dfrac{{2({{(1 - x)}^{2 - \alpha }} + {x^{2 - \alpha }})}}{{\Gamma (3 - \alpha )}}} \right]\\
\,\,\,\,\,\,\,\,\,\,\,\,\,\,\, + \dfrac{{{K_\beta }{t^\alpha }{e^{\beta t}}}}{{2\cos (\dfrac{{\beta \pi }}{2})}}\left[ {\dfrac{{24({{(1 - x)}^{4 - \beta }} + {x^{4 - \beta }})}}{{\Gamma (5 - \beta )}} - \dfrac{{12({{(1 - x)}^{3 - \beta }} + {x^{3 - \beta }})}}{{\Gamma (4 - \beta )}} + \dfrac{{2({{(1 - x)}^{2 - \beta }} + {x^{2 - \beta }})}}{{\Gamma (3 - \beta )}}}. \right]
\end{array}\]
The exact solution of Eq.(\ref{eq000002}) is given by $u(x,t) = {t^\alpha }{e^{\beta t}}{x^2}{(1 - x)^2}.$\\
\begin{table}[ht]
\centering
\caption{The absolute errors of the present scheme and
numerical methods in Yang et al.\cite{yang2010numerical} with $\alpha  = 0.4,\,\,\beta  = 1.7,\,\,\,t = 2$ for \ref{eq000002}.}
\begin{tabular}[t]{ccccccc}
\hline
$h$&Proposed Method&L1/L2\cite{yang2010numerical}&Standard/shifted Grünwald\cite{yang2010numerical}\\
\hline
$\frac{1 }{{50}}$&5.9364e-004&1.8144e-002&2.8191e-003\\
$\frac{1 }{{100}}$&1.5243e-004&9.4917e-003&1.5093e-003\\
$\frac{1 }{{200}}$&3.8468e-005&4.8584e-003&7.7821e-004\\
$\frac{1 }{{400}}$&9.6028e-006&2.4586e-003&3.9459e-004\\
\hline
\end{tabular}
\end{table}%
The absolute errors of $L1/L2$ approximation method and the standard/shifted Grünwald approximation method given in Yang et al.\cite{yang2010numerical} and presented scheme are shown in Tables 2 and they are compared for different values of $h$ and $\alpha  = 0.4,\,\,\beta  = 1.7,\,\,\,t =2$ and ${K_\alpha } = {K_\beta } = 2$.\\
{\bf Example 3.} We consider the following Riesz space fractional advection-dispersion equation with a source term and delay
\begin{equation}\label{eq80}
\dfrac{{\partial u(x,t)}}{{\partial t}} +\dfrac{{{\partial ^\gamma }u(x,t - 1)}}{{\partial {{t}^\gamma }}}= \dfrac{{{\partial ^\alpha }u(x,t)}}{{\partial {{\left| x \right|}^\alpha }}} + \dfrac{{{\partial ^\beta }u(x,t)}}{{\partial {{\left| x \right|}^\beta }}} + f(x,t),
\end{equation}
subject to the initial condition:
\begin{equation}\label{eq81}
\begin{array}{l}
u(x,t) ={x^2}{(1 - x)^2}{t^2},\,\,\, - 1 \le t \le 0,\,\,\,0 \le x \le 1,\\
u(0,t) = 0,\,\,u(1,t) =0,\,\,0 \le t \le 1,
\end{array}
\end{equation}
where $0 < \gamma  < 1,\,\,0 < \alpha  < 1,\,\,1 < \beta  \le 2$ and
\[\begin{array}{l}
f(x,t) = 2{x^2}{(1 - x)^2}t + \left[ {\dfrac{{\Gamma (3)}}{{\Gamma (3 - \gamma )}}{t^{2 - \gamma }} - 2\dfrac{{\Gamma (2)}}{{\Gamma (2 - \gamma )}}{t^{1 - \gamma }} + \dfrac{{\Gamma (1)}}{{\Gamma (1 - \gamma )}}{t^{ - \gamma }}} \right]{x^2}{(1 - x)^2}\\
\,\,\,\,\,\,\,\,\,\,\,\,\,\,\, + \dfrac{{{t^2}}}{{2\cos (\dfrac{{\alpha \pi }}{2})}}\left[ {\dfrac{{24({{(1 - x)}^{4 - \alpha }} + {x^{4 - \alpha }})}}{{\Gamma (5 - \alpha )}} - \dfrac{{12({{(1 - x)}^{3 - \alpha }} + {x^{3 - \alpha }})}}{{\Gamma (4 - \alpha )}} + \dfrac{{2({{(1 - x)}^{2 - \alpha }} + {x^{2 - \alpha }})}}{{\Gamma (3 - \alpha )}}} \right]\\
\,\,\,\,\,\,\,\,\,\,\,\,\,\,\, + \dfrac{{{t^2}}}{{2\cos (\dfrac{{\beta \pi }}{2})}}\left[ {\dfrac{{24({{(1 - x)}^{4 - \beta }} + {x^{4 - \beta }})}}{{\Gamma (5 - \beta )}} - \dfrac{{12({{(1 - x)}^{3 - \beta }} + {x^{3 - \beta }})}}{{\Gamma (4 - \beta )}} + \dfrac{{2({{(1 - x)}^{2 - \beta }} + {x^{2 - \beta }})}}{{\Gamma (3 - \beta )}}} \right].
\end{array}\]
Its exact solution is $u(x,t) = {x^2}{(1 - x)^2}{t^2}$. Some numerical results are reported at Tables 3--5.
\begin{table}[h]
\centering
\caption{The absolute errors and the convergence orders of presented method (\ref{eq62}) for example 2 with $\gamma=0.5$.}
\begin{tabular}{ccccccccccc}
\hline
   ~    &  \multicolumn{2}{c}{$\beta=1.2$}   &&  \multicolumn{2}{c}{$\beta=1.6 $} && \multicolumn{2}{c}{$\beta=1.9$}     \\
\cmidrule{2-3} \cmidrule{5-6} \cmidrule{8-9}
   ~    &$Error$& $Order$      &&$Error $& $Order $  &&$Error $& $Order $      \\
\hline
$\alpha=0.2$, $h=\kappa  $\\
1/16  & 3.8676e-004 & --  && 5.3220e-004  & --  && 6.0947e-004  & --\\
1/32  & 1.4692e-004 & 1.40  && 1.3873e-004  & 1.94  && 1.6150e-004  & 1.92\\
1/64  & 4.4836e-005 & 1.71  && 3.4679e-005  & 2.00  && 4.1138e-005  & 1.97 \\
1/128 & 1.2423e-005 & 1.85  && 8.5082e-006  & 2.02  && 1.0283e-005  & 2.00 \\
\hline
$\alpha=0.5$, $h=\kappa  $\\
1/16  & 2.4977e-004 & --  && 4.2841e-004  & --  && 5.3059e-004  & --\\
1/32  & 1.1031e-004 & 1.18  && 1.1294e-004  & 1.92  && 1.4195e-004  & 1.90\\
1/64  & 3.6830e-005 & 1.58  && 2.8289e-005  & 2.00  && 3.6317e-005  & 1.97 \\
1/128 & 1.0832e-005 & 1.77  && 7.1049e-006  & 2.00  && 9.0898e-006  & 2.00 \\
\hline
\end{tabular}
\end{table}
\begin{table}[h]
\centering
\caption{The absolute errors and the convergence orders of presented method (\ref{eq62}) for example 2 with $\beta=1.9$.}
\begin{tabular}{ccccccccccc}
\hline
   ~    &  \multicolumn{2}{c}{$\alpha=0.2$}   &&  \multicolumn{2}{c}{$\alpha=0.5 $} && \multicolumn{2}{c}{$\alpha=0.8$}     \\
\cmidrule{2-3} \cmidrule{5-6} \cmidrule{8-9}
   ~    &$Error$& $Order$      &&$Error $& $Order $  &&$Error $& $Order $      \\
\hline
$\gamma=0.2$, $h=\kappa  $\\
1/16  & 6.0947e-004 & --  && 5.3059e-004  & --  && 2.9416e-004  & --\\
1/32  & 1.6150e-004 & 1.92  && 1.4195e-004  & 1.90  && 9.3235e-005  & 1.67\\
1/64  & 4.1138e-005 & 1.97  && 3.6317e-005  & 1.97  && 2.5705e-005  & 1.86 \\
1/128 & 1.0283e-005 & 2.00  && 9.0898e-006  & 2.00  && 6.6475e-006  & 1.95 \\
\hline
\end{tabular}
\end{table}
\begin{table}[h]
\centering
\caption{The absolute errors and the convergence orders of presented method (\ref{eq62}) for example 2 with $\alpha=0.1$.}
\begin{tabular}{ccccccccccc}
\hline
   ~    &  \multicolumn{2}{c}{$\beta=1.2$}   &&  \multicolumn{2}{c}{$\beta=1.6 $} && \multicolumn{2}{c}{$\beta=1.9$}     \\
\cmidrule{2-3} \cmidrule{5-6} \cmidrule{8-9}
   ~    &$Error$& $Order$      &&$Error $& $Order $  &&$Error $& $Order $      \\
\hline
$\gamma=0.7$, $h=\kappa  $\\
1/16  & 4.0456e-004 & --  && 5.5129e-004  & --  && 6.2364e-004  & --\\
1/32  & 1.5142e-004 & 1.42  && 1.4380e-004  & 1.94  && 1.6524e-004  & 1.92\\
1/64  & 4.5762e-005 & 1.73  && 3.5984e-005  & 2.00  && 4.2094e-005  & 1.97 \\
1/128 & 1.2597e-005 & 1.86  && 8.8387e-006  & 2.02  && 1.0524e-005  & 2.00 \\
\hline
\end{tabular}
\end{table}
\section{Conclusion}
In this paper, the fractional backward differential formulas method and shifted Gr\"{u}nwald difference (WSGD) operators are used to approximate the Riesz space fractional derivative.
The presented numerical method is proposed and applied to solve the Riesz space fractional advection-dispersion equations with delay. We prove that the difference scheme is conditionally stable and convergent with the accuracy ${\rm O}({\kappa ^2} + {h^2})$. Furthermore, the analytical solution for RFADED in terms of the functions of Mittag--Leffler type has been obtained. Finally, some results are given to demonstrate the effectiveness of theoretical analysis. In the future, we would like to investigate this method and technique of high order for solving this problem.
\section*{References}

\bibliography{mybibfile}

\end{document}